\documentclass[preprint]{imsart}
\usepackage{amsthm,amsmath}
\usepackage{amssymb,euscript}
\usepackage{xcolor}
\RequirePackage[numbers]{natbib}
\RequirePackage[colorlinks,citecolor=blue,urlcolor=blue]{hyperref}
\usepackage{graphicx}
\usepackage{pgfplots}
\usepgfplotslibrary{fillbetween}
\usepackage{csquotes}
\usepackage[caption=false]{subfig}
\pgfplotsset{compat=1.10}

\doi{}
\pubyear{0000}
\volume{0}
\firstpage{1}
\lastpage{30}

\startlocaldefs

\newtheorem{lemma}{Lemma}
\newtheorem{theorem}{Theorem}
\newtheorem{definition}{Definition}
\newtheorem{corollary}{Corollary}
\newtheorem{assumption}{Assumption}

\newtheorem{Example}{Example}
\def\a{\alpha}

\def\de{\delta}
\def\De{\Delta}

\def\ga{\gamma}

\def\la{\lambda}

\def\ep{\varepsilon}

\def\cov{{\rm Cov}}

\def\sp{{\rm supp}}
\def\var{{\rm Var}}

\def\In{{\rm 1 \hskip-2.9truept l}}

\newcommand{\wh}{\widehat}

\def\N{{\mathbb N}}
\def\R{{\mathbb R}}

\def\eps{\varepsilon}

\def\E{{\mathbb E}}
\def\I{{\cal I}}

\def\Z{{\mathbb Z}}

\def\V{{\cal V}}
\def\ovl{\overline}
\def\wh{\widehat}

\newenvironment{proof1}{%
    \vspace{0.3cm} \pagebreak [2]%
    \par%
    \noindent%
    {\it Proof~}}{\qed}%
\newenvironment{remark}{%
    \vspace{0.3cm} \pagebreak [2]%
    \par%
    \refstepcounter{proposition}
    \noindent%
    {\bf Remark~\theproposition.\ }}{ }%
\endlocaldefs
\begin{document}

\begin{frontmatter}

\title{Asymptotic normality of simultaneous estimators of cyclic long-memory processes}
\runtitle{CLT for cyclic long-memory processes}


\begin{aug}
\author{\fnms{Antoine} \snm{Ayache}
\ead[label=e1]{Antoine.Ayache@univ-lille.fr}}
\address{Laboratoire Paul-Painlev\'e (UMR CNRS 8524), Universit\'e de Lille, B\^atiment M2, \\ Cit\'e Scientifique, 59655 Villeneuve d'Ascq, France \\
\printead{e1}}
\author{\fnms{Myriam} \snm{Fradon}\ead[label=e2]{Myriam.Fradon@univ-lille.fr}}
\address{Laboratoire Paul-Painlev\'e (UMR CNRS 8524), Universit\'e de Lille, B\^atiment M2, \\ Cit\'e Scientifique, 59655 Villeneuve d'Ascq, France \\
\printead{e2}}
\author{\fnms{Ravindi} \snm{Nanayakkara}\ead[label=e3]{D.Nanayakkara@latrobe.edu.au}}
\address{Department of Mathematics and Statistics, La Trobe University, Melbourne, 3086, Australia. \\
\printead{e3}}
\author{\fnms{Andriy} \snm{Olenko}
\thanksref{t1} \ead[label=e4]{A.Olenko@latrobe.edu.au}}
\thankstext{t1}{Corresponding author.}
\address{Department of Mathematics and Statistics, La Trobe University, Melbourne, 3086, Australia. \\
\printead{e4}}

\runauthor{Ayache et al.}

\end{aug}

\begin{abstract}
Spectral singularities at non-zero frequencies play an important role in investigating cyclic or seasonal time series. The publication \cite{AAFO} introduced the generalized filtered method-of-moments approach to simultaneously estimate singularity location and long-memory parameters.  This paper continues studies of these simultaneous estimators. A wide class of Gegenbauer-type semi-parametric models is considered.  Asymptotic normality of several statistics of the cyclic and long-memory parameters is proved. New adjusted estimates are proposed and investigated. The theoretical findings are illustrated by numerical results. The methodology includes wavelet transformations as a particular case.
\end{abstract}


\begin{keyword}
\kwd{Central limit theorem}
\kwd{cyclic long-memory}
\kwd{filter}
\kwd{wavelet}
\kwd{estimators of parameters}
\kwd{asymptotic normality}
\end{keyword}


\end{frontmatter}

 \section{Introduction} 

Time series with cyclic long-memory behaviours attracted increasing attention in recent years, see~\cite{AAFO, Arteche:2020, ArtRob:1999, ArtRob:2000, Castro:2020} and the references therein. It was due to importance of such time series in finance, hydrology, cosmology, internet modelling, and other applications to data with non-seasonal cyclicities, see~\cite{Arteche:2020, ArtRob:1999, Arteche:2011, Boubaker:2015, Ferrara:2001, Whitcher:2004}. At the same time, various statistics of cyclic long-memory processes have complex asymptotic behaviour that has not yet been fully understood and investigated, see~\cite{Hosoya:1997, Ivanov:2013, Klyolen:2012, Olenko:2013}.

To link characterizations of the long-memory phenomena in temporal and spectral domains researchers usually employ Abelian and Tauberian theorems. These results establish connections between asymptotics of covariance functions at the infinity and singularities of the corresponding spectral densities, see~\cite{Klyolen:2012, LeonOle:2013}. The most frequent definition of long-memory in the literature is a hyperbolic-type decay of a non-integrable covariance function. While this classical long-memory dependence is often related to unboundedness of spectral densities at the origin, spectral singularities at nonzero frequencies can also result in hyperbolic-type oscillating non-integrable covariance functions. Such spectral representations can be used to simultaneously model cyclicity and long-memory. 

Cyclical long-memory time series are much more difficult to investigate and there were relatively few publications on this topic compared to classical models with the only singularity at the~origin. Several least squares and likelihood-based approaches have been proposed to estimate parameters of singularity poles, see~\cite{Arteche:2020, ArtRob:1999, ArtRob:2000, Barboza:2017, Beran:2009, Espejo:2015, Giraitis:2001, Hidalgo:2005, Tsai:2015}. Unfortunately, for the majority of these approaches incorrect specifications of a statistical model can result in inconsistent estimates of the parameters. The empirical studies in \cite{RePEc:pra:mprapa:96313, Whitcher:2004} demonstrated various issues of the traditional estimators and that wavelet-based approach can give results that are equivalent to ordinary least squares and maximum likelihood estimates under the assumption of knowing the explicit form of the spectrum. However, for the cases when the model is not fully specified, wavelets can provide better estimates.

To avoid repetitions, we refer the readers to  very detailed motivation, discussion and various examples in \cite{AAFO}.

This paper investigates time series which spectral density $f(\cdot)$ has the following semiparametric form
\[ f(\lambda) = \frac{h(\lambda)}{{|{\lambda}^2 - s_0^2|}^{2\alpha}}, \quad \lambda \in \mathbb{R}. \]

The parameter $s_0$ determines cyclic behaviour while $\alpha$ is a long-memory parameter. For example, the Gegenbauer model \cite{Espejo:2015} has a spectral density of this form.

\begin{figure}[!htb]
    \centering  \vspace{-0.3cm}
    \subfloat[Realization]{\label{fig21a}
    \includegraphics[trim={0cm 0cm 0cm 0cm},clip,width=0.45\textwidth, height=0.28\textheight]{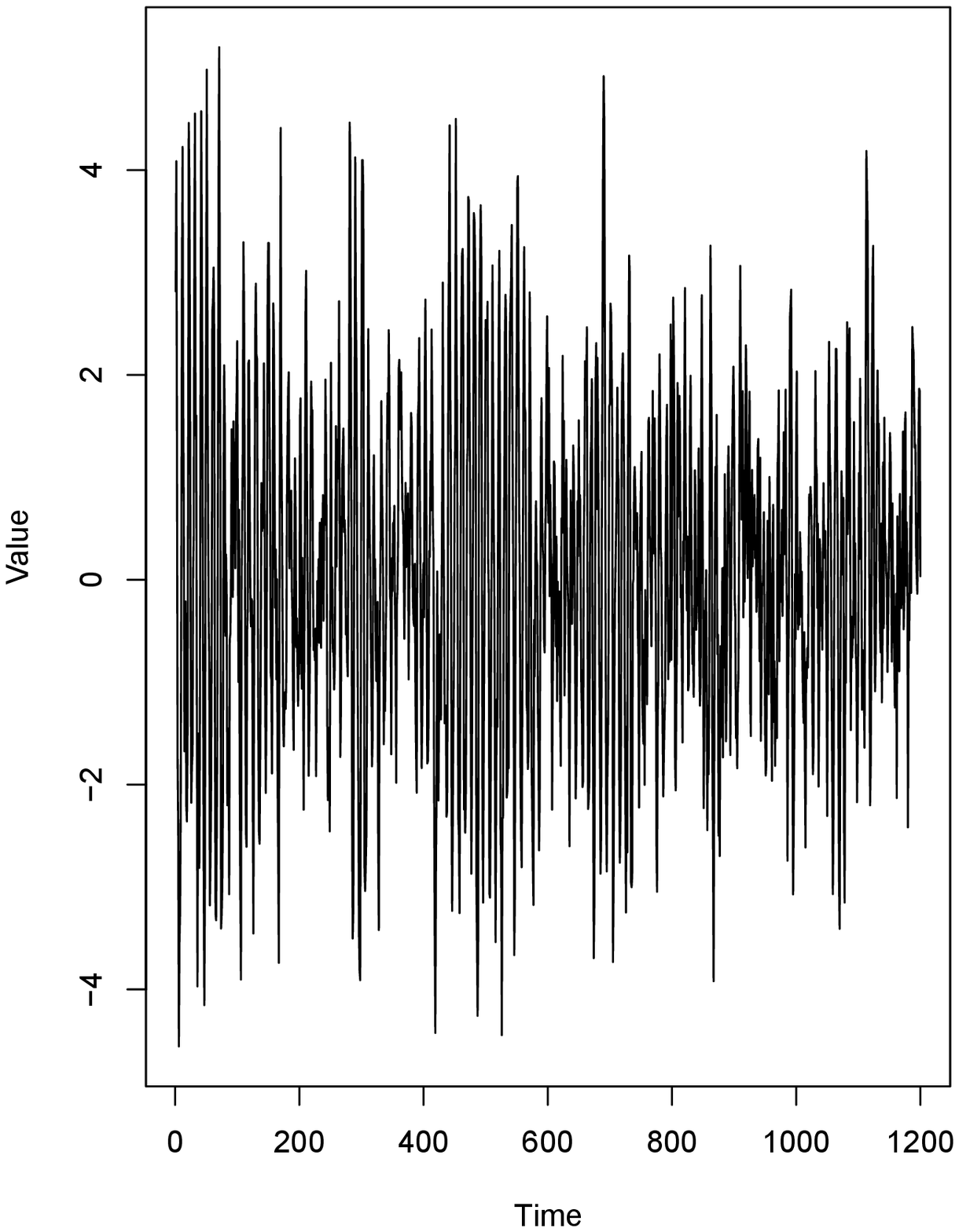}}
    \centering
    \subfloat[Periodogram]{\label{fig21b}
    \includegraphics[trim={0cm 0cm 0cm 0cm},clip,width=0.45\textwidth, height=0.28\textheight]{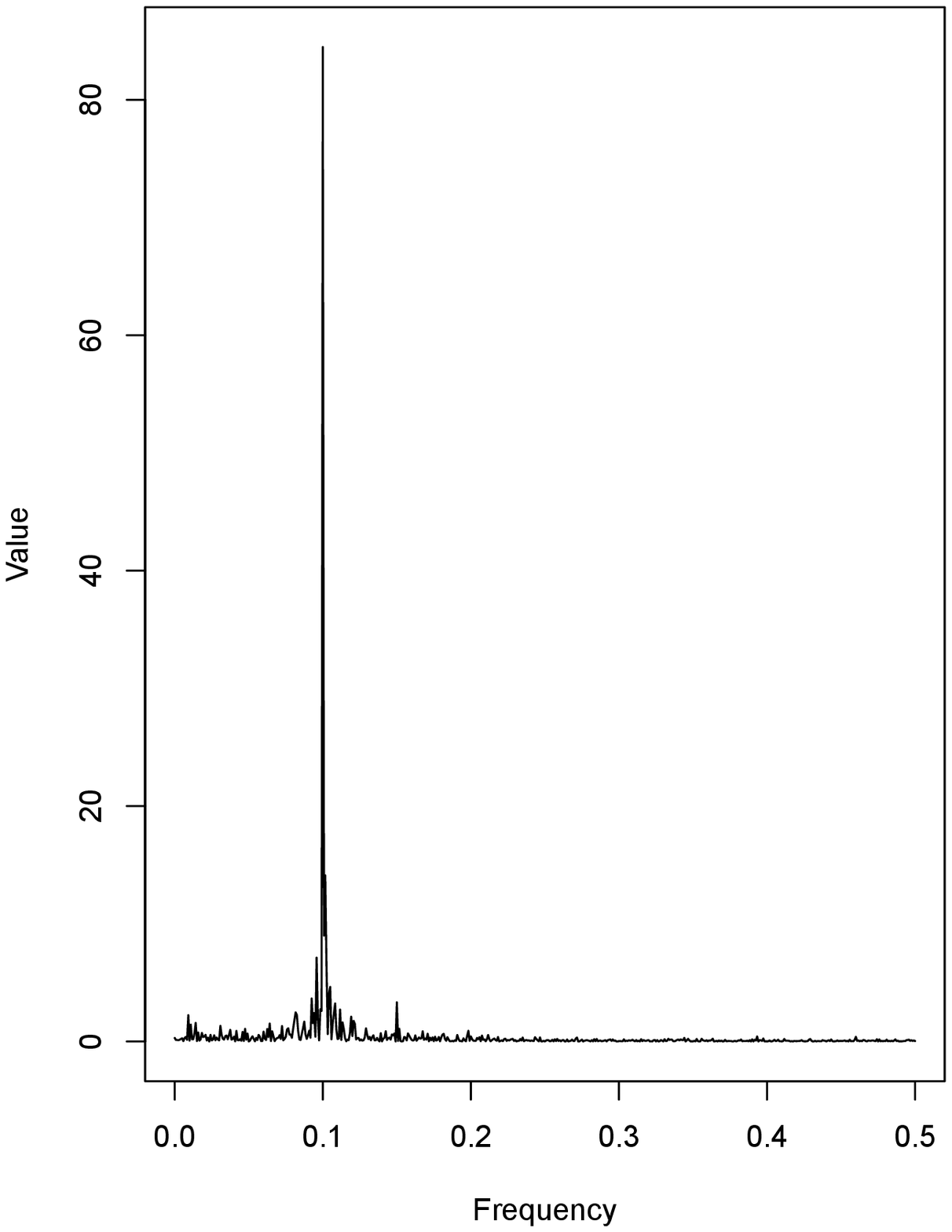}}\\
    \centering \vspace{-0.4cm}
    \subfloat[Sample covariance function]{\label{fig21c}
    \includegraphics[trim={0cm 0cm 0cm 0cm},clip,width=0.45\textwidth, height=0.28\textheight]{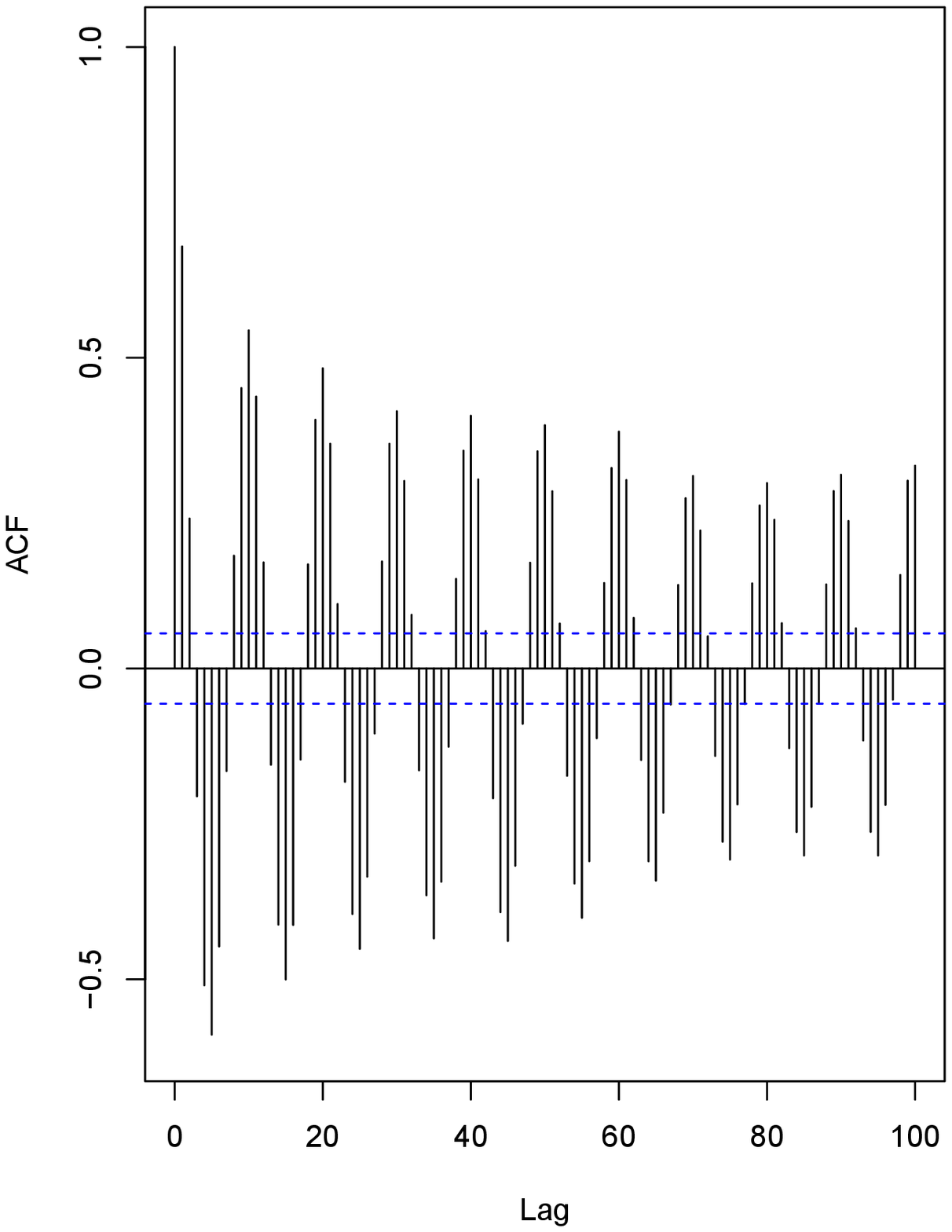}}
    \centering 
    \subfloat[Wavelet coefficients]{\label{fig21d}
    \includegraphics[trim={0.2cm 0.6cm 1.1cm 1.1cm},clip,width=0.45\textwidth, height=0.28\textheight]{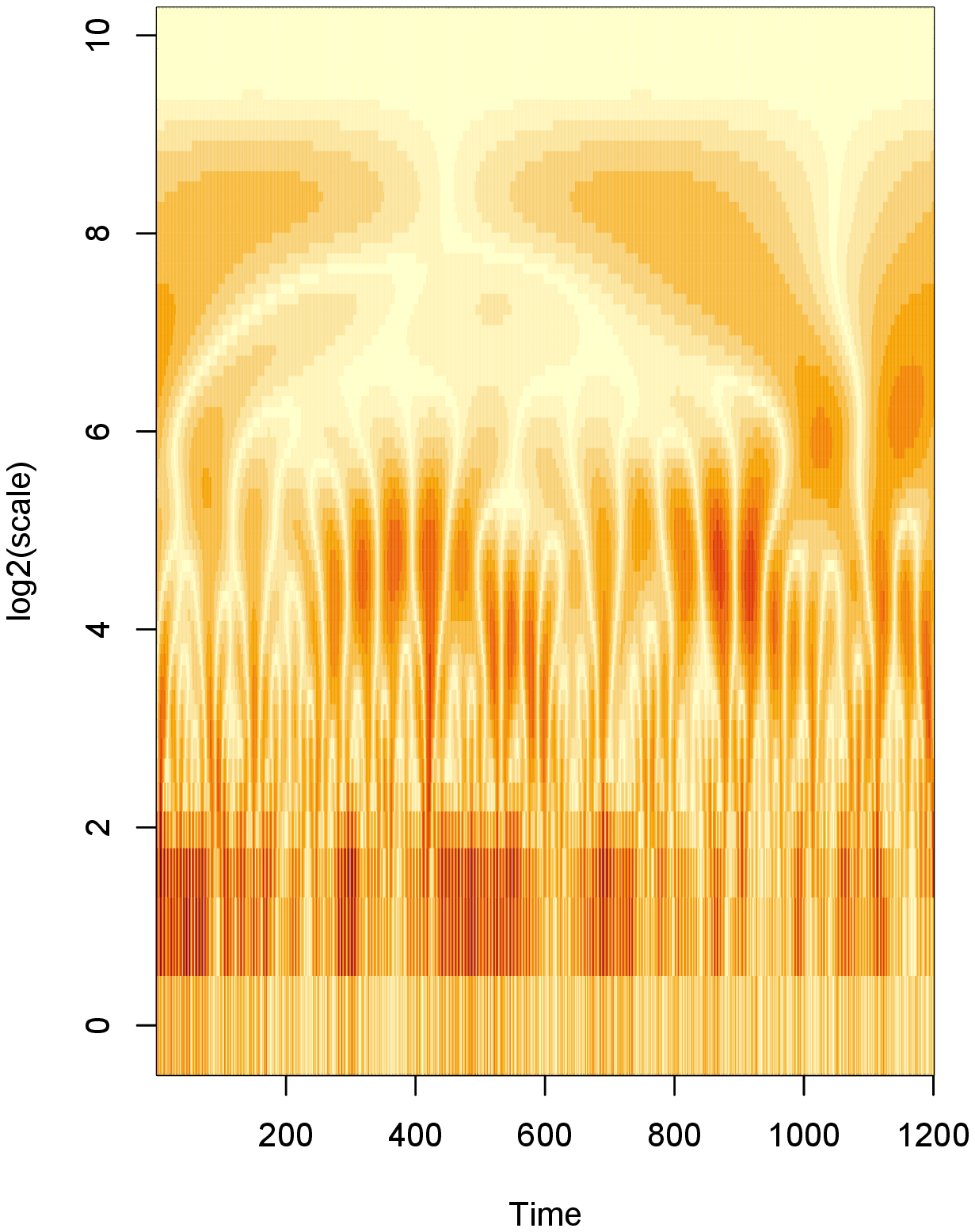}}
    \caption{Cyclic long-memory time series}
    \label{fig21}
\end{figure}

Figure~\ref{fig21} shows a realization of such time series together with its estimated spectral density and  covariance function. In this example a spectral density with a sharp spike at its singularity location was chosen.  It clear demonstrates that the spectral density has a singularity at a non-zero frequency and the corresponding covariance function indicates some cyclic behaviour. The wavelet coefficients of this time series are shown in the fourth subplot. Unfortunately, contrary to perfect cyclic signals or spectral densities with singularity at the origin, it is more difficult to use the wavelet approach for estimating cyclicity  and long-memory parameters simultaneously. An even more challenging problem is a development of statistical inference for  these parameters.

The publication \cite{AAFO} proposed a new methodology for simultaneous estimation of cyclic and long-memory parameters. It used filter transformations of functional time series. The approach included wavelet transformations as a particular case. The consistency of the proposed estimators was proved.

This paper further develops the approach from \cite{AAFO}. Now we obtain asymptotic normality of the proposed estimators. It requires very careful investigations of quadratic functionals of filter coefficients and their increments. Obtaining asymptotic properties of wavelet-based statistics is a difficult problem and there are only few general results about their asymptotic normality. The developed methodology and the obtained results can also find applications for other wavelet-based statistics.

In addition, for the case when empirical values of the statistics are outside the feasible region, we propose  new adjusted estimators and investigate their properties. It is shown that these estimators have same asymptotic distributions as the corresponding ones in \cite{AAFO}, but are computationally simpler.

The article is organized as follows. Section~\ref{sec_2} gives basic definitions and introduces a semi-parametric model and filter transforms studied in this paper. Various asymptotic properties of quadratic functionals of filter transforms are derived in Section~\ref{sec_3}. 
Section~\ref{sec_4} proves asymptotic normality of two auxiliary statistics of the semiparametric model, which are based on quadratic functionals of filter transforms and their increments. Section~\ref{sec_5} proposes and investigates adjusted simultaneous estimators of the location and long-memory parameters. Numerical studies to support the theoretical findings are presented in Section~\ref{sec_6}.

All computations, plotting and simulations in this article were performed using the software~R version~4.0.3 and Maple 17, Maplesoft. In particular, the R packages~{\sc waveslim} \cite{Waveslim:2020} and {\sc MassSpecWavelet} \cite{Massspecwavelet:2006} were used to simulate realizations of cyclic long-memory processes and compute their wavelet  transforms in  the numerical examples. A reproducible version of the code in this paper is available in the folder \enquote{Research materials} from the  website \url{https://sites.google.com/site/olenkoandriy/}.
 \section{Definitions and assumptions}\label{sec_2}

This section introduces classes of functional time series and their filter transforms that are used in the paper. The notations are consistent with ones in~\cite{AAFO}, where the authors proposed simultaneous filter estimators of parameters of cyclic long-memory processes.

In the following $\{a_j\}_{j\in \N}$ denotes an arbitrary unboundedly strictly monotone increasing sequence of positive real numbers.  $\{m_j\}_{j\in\N}$ is an unboundedly increasing sequence of positive integers.  $\{b_{jk}\}_{(j,k)\in \N\times \Z}$ stands for an infinite array of real numbers. 

The symbols  $\xrightarrow {a.s.}$ and $\xrightarrow {d}$ will be used for almost sure convergence and convergence in distribution respectively.

Let $ X(t),$ $t\in \mathbb{R},$ be a measurable mean-square continuous real-valued stationary zero-mean 
Gaussian stochastic process 
on a probability space $(\Omega, \mathcal{F}, P),$ with the 
covariance function
\[B(r):=\cov(X(t), X(t'))=\int_ \mathbb{R} \mathrm{e}^{iu(t-t')} F(du), \quad t, t' \in \mathbb{R}, \] 
where $r=t-t'$ and $F(\cdot)$ is a non-negative finite measure on $\mathbb{R}.$
\begin{definition}
The random process $X(t),$ $t \in \mathbb{R},$ possesses an absolutely
 continuous spectrum if there exists a non-negative function $f(\cdot)\in L_{1}( \mathbb{R})$ such that
\[F(u)=\int_{-\infty }^ {u} f(\lambda) d\lambda, \quad  u\in \mathbb{R}. \]
\end{definition}
The function $f(\cdot)$ is called the spectral density of the process $ X(t).$

The process $ X(t), t \in \mathbb{R},$ with an absolutely continuous spectrum has the 
following isonormal spectral representation
\[X(t)= \int_ {\mathbb{R}}\mathrm{e}^{i t\lambda} \sqrt{f(\lambda)} dW(\lambda),\]
where $W(\cdot)$ is a complex-valued Gaussian orthogonal random measure on $\mathbb{R}.$

For a real-valued process $X(t)$ the function $f(\cdot)$ is  even and the random measure  $W(\cdot)$ satisfies the condition 
 $W\left(\left[\lambda_1,\lambda_2\right]\right)=W\left(\left[-\lambda_2,-\lambda_1\right]\right)$ for any $\lambda_2>\lambda_1>0,$ see \cite[\S 6]{Taqqu:1979}.
 
 The following assumption in the spectral domain introduces the semi-para\-met\-ric model investigated in this paper. 
\begin{assumption}\label{Assumption_1}
Let the spectral density $f(\cdot)$ of $X(t)$ admit the following representation
\[f(\lambda)=\frac{h(\lambda)}{|\lambda^{2}-s_{0}^{2}|^{2\alpha}},\quad \lambda \in \mathbb{R},\]
where $s_{0}> 1,\, \alpha\in (0,{1}/{2})$ and $h(\cdot)$ is an even non-negative bounded function that is  four times continuously differentiable. Its derivatives of order $i$ satisfy $h^{(i)}(0)=0,$ $i=1,2,3,4.$ Also, $h(0)=1,$ $h(\cdot)>0$ in some neighborhood of $\lambda=\pm  s_0,$ and for all $\varepsilon>0$ it holds 
\[\int_\mathbb{R}\frac{h(\lambda)}{(1+|\lambda|)^\varepsilon }d\lambda<\infty.\]
\end{assumption}

Stochastic processes with spectral densities satisfying Assumption $\ref{Assumption_1}$ exhibit cyclic long memory. The boundedness of $h(\cdot)$ guarantees that their spectral densities have singularities only at the locations $ \pm s_{0}.$ Covariance functions of 
such processes are unintegrable and have hyperbolically decaying oscillations when  $\alpha\in\left(0,1/2\right),$ see \rm\cite{ArtRob:1999}.
For example, the Gegenbauer random processes satisfy Assumption~\rm\ref{Assumption_1}, see \rm\cite{Espejo:2015}.

Real-valued functions $\psi (t)\in L_1(\mathbb{R}),$ $t \in \mathbb{R},$ are used to introduce filter transforms of  the process $X(t)$. The Fourier transform  $\widehat{\psi}$ is  defined, for each $\lambda \in \mathbb{R}$, as $\widehat{\psi}(\lambda)=\int_\mathbb{R}e^{-i\lambda t} \psi(t) dt.$
It follows from properties of $\psi(\cdot) $ that $\widehat{\psi}(\cdot)$ is a bounded even function.

\begin{assumption}\label{Assumption_2}
Let $ {\rm supp} \, \widehat \psi\subset [-A,A],$ $A>0,$ and
 $\widehat \psi(\cdot)$ is of bounded variation on~$[-A,A].$ 
 \end{assumption}
 This assumption is technical and can be replaced by a sufficiently fast decay rate of $\widehat \psi(\cdot)$ at infinity.
 
\begin{definition}\label{def2} The filter transform of the process $X(t)$ is the array of  centred real-valued Gaussian random variables
$\{\de_{jk}\}_{(j,k)\in \N\times \Z}$  defined as
\begin{equation}
\label{eq:de}
\de_{jk}:=\frac{1}{\sqrt{a_j}} \int_{\mathbb{R}} \psi\bigg(\frac{t-b_{jk}}{a_j}\bigg)X(t) dt=\sqrt{a_j}\int_\R e^{i b_{jk}\xi}\,\frac{\ovl{\wh{\psi}(a_j\xi)}\sqrt{h(\xi)}}{|\xi^2-s_0^2|^\a}\, dW(\xi).
\end{equation}
\end{definition}
Definition~\ref{def2} provides equivalent expressions of the filter transform in the spectral and time domains.

It is easy to see that
\begin{equation}
\label{varde}
\var (\de_{jk})  =a_j\int_\R \frac{\big |\wh{\psi}(a_j\xi)\big |^2 h(\xi)}{|\xi^2-s_0^2|^{2\a}}\, d\xi.
\end{equation}

To guarantee that at each level $j \in \N$ the sequence $\{b_{jk}\}_{k \in \Z}$ does not have concentration points and covers all spectral range the following  assumption is rather standard in the literature.
\begin{assumption}\label{Assumption_3}
For all $j\in\N$ and for every $(k,l)\in\Z^2$ it holds 
\begin{equation}
\label{eq:bgam}
|b_{jk}-b_{jl}|\ge \ga_j |k-l|,
\end{equation}
where $\{\ga_j\}_{j\in\N}$ is a sequence of positive real numbers.
\end{assumption}
To get exact asymptotic behaviours of the considered statistics few versions of this assumption will be more precisely specified later.

 A very detailed motivation, discussion,  and various particular examples, that include wavelet transforms and Gegenbauer processes as special important cases, can be found in~\cite{AAFO}.

\section{Preliminary results}\label{sec_3}

This section derives some properties of the filter transforms and their variances that will be used in the following sections to obtain the CLT for simultaneous estimators of cyclic long-memory parameters. 

Let 
\begin{equation}
\label{prop:ctl1:eq1}
\de_{j}^{(2,m_j)}:=\sum_{k=1}^{m_j} \de_{jk}^2, \quad j\in\N.
\end{equation}

\begin{theorem}
\label{prop:ctl1}
Assume that 
\begin{equation}
\label{prop:ctl1:eq2}
\lim_{j\rightarrow +\infty} \frac{a_j \log (m_j)}{\ga_j\,m_j ^{1/2}}=0.
\end{equation}
Then, when $j\rightarrow +\infty$, the random variables
\begin{equation}
\label{prop:ctl1:eq3}
Y_j:=\frac{\de_{j}^{(2,m_j)}-\E(\de_{j}^{(2,m_j)})}{\sqrt{\var (\de_{j}^{(2,m_j)})}}
\end{equation}
converge in distribution to a standard Gaussian random variable.
\end{theorem}

To derive Theorem~\ref{prop:ctl1} we will use the following three lemmas. The first lemma is obtained by applying the Taylor-Lagrange formula, the second one is a rather known result and the third statement was proved in \cite{AAFO}.

Let the function $\I_\zeta(\cdot),$ $\zeta\in\R,$ be defined for  $x\in \big [-(2A)^{-1},(2A)^{-1}\big ]$ as
\begin{equation}
\label{lem:var-de:eq1}
\I_\zeta (x) := \int_\R e^{i\zeta\eta} \frac{|\wh{\psi}(\eta)|^2 h(x\eta)}{\big(s_0^2-x^2\eta^2\big)^{2\a}} \, d\eta.
 \end{equation}
\begin{lemma}
\label{lem:var-de}
If Assumptions~{\rm\ref{Assumption_1}} and {\rm \ref{Assumption_2}} hold true, then $\I_\zeta(x)$ is four times continuously differentiable with respect to $x,$ and there is a finite constant $c_1>0$ (not depending on $\zeta$ and $x$) such that, for all $\zeta\in\R$ and $|x|\le(2A)^{-1},$ it holds
\begin{equation}
\label{lem:var-de:eq2}
\bigg | \I_\zeta (x)-s_0^{-4\a}\int_{\R}e^{i\zeta\eta}|\wh{\psi}(\eta)|^2\,d\eta-2\a s_0^{-4\a-2}\int_{\R}e^{i\zeta\eta}\eta^2|\wh{\psi}(\eta)|^2\,d\eta\cdot x^2\bigg|\le c_1 \,x^4.
\end{equation}
\end{lemma}

\begin{proof1}{\it of Lemma {\rm \ref{lem:var-de}.}}  
Note that $\I_\zeta(\cdot)$ is a real-valued function since $\wh{\psi}(\cdot)$ and $h(\cdot)$ are even real-valued functions. It follows from (\ref{lem:var-de:eq1}), Assumptions~\ref{Assumption_1} and \ref{Assumption_2}  that 
 \[
\I_\zeta (x)  =\int_{-A}^A e^{i\zeta\eta}\frac{|\wh{\psi}(\eta)|^2 h(x\eta)}{\big(s_0^2-x^2\eta^2\big)^{2\a}} \, d\eta =\int_{-A}^A e^{i\zeta\eta} |\wh{\psi}(\eta)|^2 f(\eta x) \, d\eta.
\]

To use the Taylor formula for $\I_\zeta(x)$ when $ x\in \left[-(2A)^{-1},(2A)^{-1}\right]$ one notes that
$x\in \left[-(2A)^{-1},(2A)^{-1}\right]$ and $\eta\in[-A,A]$ imply $|\eta x|\le 1/2$ and $s_0^2 - \eta^2 x^2 > 3/4$ since $s_0>1$. 
As by Assumption~\ref{Assumption_1} the function $h(\cdot)$ is four times continuously differentiable, hence $f(\cdot)$ has four continuous derivatives with respect to $x$ on $\left[-(2A)^{-1},(2A)^{-1}\right]$ for any fixed $\eta$ in $[-A,A]$. To prove that $I_\zeta(\cdot)$ is four times continuously differentiable, it is enough to show that the corresponding integrand and its first four derivatives with respect to $x$ are dominated by integrable functions that do not depend on $x$. 

First, for the integrand in (\ref{lem:var-de:eq1}) we get
\[
\left| e^{i\zeta\eta} |\wh{\psi}(\eta)|^2 f(\eta x) \right| \le \left( \frac43 \right)^{2\a} |\wh{\psi}(\eta)|^2 \sup_{y\in [-1/2,1/2]} |h(y)|,\]
where the right hand side is bounded and therefore integrable on  $[-A,A].$

The $n^\text{th}$ derivative of the function $f(\eta x)$ with respect to $x$ satisfies
\[
\left| \frac{\partial^n }{\partial x^n} f(\eta x)\right|
= \left| \sum_{k=0}^n  \binom{n}{k} ~\eta^{n-k} ~h^{(n-k)}(x\eta) 
         ~\frac{\partial^k}{\partial x^k} \left( (s_0^2 - \eta^2 x^2)^{-2\a} \right) \right| \]
         \[
\le \sum_{k=0}^n \binom{n}{k} ~A^{n-k} ~\sup_{y\in [-1/2,1/2]} \big | h^{(n-k)}(y)\big | 
         \left| \frac{\partial^k}{\partial x^k} \left( (s_0^2 - \eta^2 x^2)^{-2\a} \right) \right|.        
\]
For $k$ in $\{1,2,3,4\}$ we provide very simple convenient bounds for the derivatives in the last expression, which will be useful later:
\begin{equation} \label{bound1}
\left| \frac{\partial}{\partial x} \left( (s_0^2 - \eta^2 x^2)^{-2\a} \right) \right|
 = \left| \frac{ 4\a\eta^2 x}{ (s_0^2 - \eta^2 x^2)^{2\a+1}} \right| 
 \le 4\a \frac{A}{2} \left( \frac43 \right)^{2\a+1} \le 2A,
\end{equation}
\[
\left| \frac{\partial^2}{\partial x^2} \left( (s_0^2 - \eta^2 x^2)^{-2\a} \right) \right|
= \left| 4 \a \eta^2 ~\frac{ (4\a+1)\eta^2 x^2 + s_0^2 }{ (s_0^2 - \eta^2 x^2)^{2\a+2} } \right| \]
\begin{equation} \label{bound2}
\le 4\a A^2 \left( \frac43 \right)^{2\a+2} \left( \frac{4\a+1}{4} + s_0^2 \right)
\le 10 A^2 s_0^2,
\end{equation}
\[
\left| \frac{\partial^3}{\partial x^3} \left( (s_0^2 - \eta^2 x^2)^{-2\a} \right) \right|
= \left| 8\a (2\a+1) \eta^4 x ~\frac{ (4\a+1) \eta^2 x^2 + 3 s_0^2 }{ (s_0^2 - \eta^2 x^2)^{2\a+3} } \right|  \le 4\a(2\a+1) A^3
\] 
\begin{equation} \label{bound3}
\times \left( \frac43 \right)^{2\a+3} \left( \frac{4\a+1}{4} + 3 s_0^2 \right)
\le 4 A^3 \left( \frac43 \right)^{4} \left( \frac{3}{4} + 3 s_0^2 \right) \le 48 A^3 s_0^2,
\end{equation}
\[
\left| \frac{\partial^4}{\partial x^4} \left( (s_0^2 - \eta^2 x^2)^{-2\a} \right) \right|
= \left| \frac{ (16\a(\a+1)+3) \eta^4 x^4 +6(4\a+3) s_0^2 \eta^2 x^2 +3 s_0^4 }{ (s_0^2-\eta^2 x^2)^{2\a+4} }\right. \]
\[ \times 8\a(2\a+1)\eta^4  \Bigg| \le 8\a(2\a+1) A^4 \left( \frac43 \right)^{2\a+4} \left( \frac{16\a(\a+1)+3}{16} +\frac{6(4\a+3)}{4} s_0^2 +3 s_0^4 \right) \]
\begin{equation} \label{bound4}
\le 8 A^4 \left( \frac43 \right)^5 \left( \frac{15}{16} +\frac{15}{2} s_0^2 +3 s_0^4 \right) 
\le 400 A^4 s_0^4.
\end{equation}

Therefore the function in the integral defining $\I_\zeta(\cdot)$ and its first four derivatives are dominated by an integrable function ($|\wh{\psi}|^2$ multiplied by a large enough constant). Thus $\I_\zeta(\cdot)$ is $\mathcal{C}^4\left(\left[-(2A)^{-1},(2A)^{-1}\right]\right)$ and its derivatives can be computed by differentiation under the integral sign. For $n$ in $\{1,2,3,4\}$ it holds
\begin{equation} \label{nthder}
\frac{d^n }{d x^n}\I_\zeta(x)
= \sum_{k=0}^n \binom{n}{k} \int_{-A}^A e^{i\zeta\eta} |\wh{\psi}(\eta)|^2 ~\eta^{n-k} ~h^{(n-k)}(x\eta) 
                                           ~\frac{\partial^k}{\partial x^k} \left( (s_0^2 - \eta^2 x^2)^{-2\a} \right) \, d\eta
\end{equation}
and the Taylor-Lagrange expansion provides
\begin{equation}
\label{lem:var-de:taylor}
\bigg| \I_\zeta(x)-\I_\zeta(0) -\I'_\zeta(0) x -\I^{''}_\zeta(0)\frac{x^2}{2!} -\I^{(3)}_\zeta(0)\frac{x^3}{3!} \bigg|
\le \sup_{y\in [-(2A)^{-1},(2A)^{-1}]} |\I^{(4)}_\zeta (y)| ~~ \frac{x^4}{4!},
\end{equation}
where
\[ 
\I_\zeta(0)  =\frac{1}{s_0^{4\a}} \int_{-A}^A e^{i\zeta\eta} |\wh{\psi}(\eta)|^2 \, d\eta,
\]
since $h(0)=1.$

By Assumptions~{\rm\ref{Assumption_1}} the derivatives $h^{(l)}(0)=0$ for $l\in\{1,2,3, 4\},$ thus
\[
\frac{d^n}{d x^n} \I_\zeta(0)
= \left.\int_{-A}^A e^{i\zeta\eta} |\wh{\psi}(\eta)|^2 ~\frac{\partial^n }{\partial x^n}\left( (s_0^2-\eta^2 x^2)^{-2\a} \right)\right\vert_{x=0} \, d\eta.
\]
By (\ref{bound1}) and (\ref{bound3}) for $n=1$ and $n=3$ the derivatives $\frac{\partial^n}{\partial x^n} \left( (s_0^2-\eta^2 x^2)^{-2\a} \right)$  vanish at $x=0$. Moreover, the expression for the second derivative in the estimate (\ref{bound2}) gives
\[
\frac{d^2 }{d x^2}\I_\zeta(0) = \frac{ 4 \a }{ s_0^{4\a+2} } \int_{-A}^A e^{i\zeta\eta} |\wh{\psi}(\eta)|^2  \eta^2  \, d\eta.
\]
It follows from the estimates (\ref{bound1})-(\ref{bound4}) that for each $k=0,...,4$ the derivative  $|\frac{\partial^k}{\partial x^k}(s_0^2 - \eta^2 x^2)^{-2\a}|$ is bounded by $400 A^k s_0^4.$   Hence, by (\ref{nthder}),  for all $x\in [-(2A)^{-1}, (2A)^{-1}]$ 
\[
\left| \frac{d^4 }{d x^4}\I_\zeta(x) \right|
\le \sup_{ \substack{y\in [-1/2,1/2]\\ n\in \{0,\ldots,4\}}} |h^{(n)}(y)| 
  ~\sum_{k=0}^4 \binom{4}{k} \int_{-A}^A |\wh{\psi}(\eta)|^2  A^{n-k}(400 A^k s_0^4) \, d\eta \]
  \[
\le 6400\, s_0^4\, A^4  \sup_{\substack{y\in [-1/2,1/2]\\ n\in \{0,\ldots,4\}}} |h^{(n)}(y)| ~\int_{-A}^A |\wh{\psi}(\eta)|^2 \, d\eta =:c_2.
\]
Finally, the estimate (\ref{lem:var-de:taylor}) becomes
\[
\bigg| \I_\zeta(x) -\frac{1}{s_0^{4\a}} \int_{-A}^A e^{i\zeta\eta} |\wh{\psi}(\eta)|^2 \, d\eta
                   -\frac{ 4\a }{ s_0^{4\a+2} } \int_{-A}^A e^{i\zeta\eta} |\wh{\psi}(\eta)|^2 ~\eta^2 \, d\eta\cdot \frac{x^2}{2!} \bigg|
\le  \frac{c_2 }{4!} \cdot x^4,
\]
which completes the proof.
\end{proof1}

The following lemma is an immediate corollary of the Gershgorin circle theorem.

\begin{lemma}
\label{lem:sp-radius}
Let $U=(u_{ij})_{1\le i, j\le n}$ be a square matrix of order $n$ with complex elements. If $\rho (U)$ is the spectral radius of $U$, that is 
\[
\rho(U):=\max\big\{|\la|:\,  \la\ \mbox{is an eigenvalue of $U$}\big\},
\]
then
\[
\rho (U)\le \min\Big\{\max_{1\le i \le n} \sum_{j=1}^n |u_{ij}|\,, \max_{1\le j \le n} \sum_{i=1}^n |u_{ij}|\Big\}.
\]
\end{lemma}

\begin{lemma}{\rm \cite{AAFO}} 
 Let Assumptions~{\rm\ref{Assumption_1}} hold true. Then there exists a finite constant $c_3$ such that, for every $j\in\N$ such that $a_j\ge 2A$   and for all $(k,l)\in\Z^2$, one has
\begin{equation}
\label{lem:cov-b:eq1}
\big |\cov(\de_{jk},\de_{jl})\big |\le c_3\Big (\In_{\{k=l\}}+\In_{\{k\ne l\}}\, a_j |b_{jk}-b_{jl}|^{-1}\Big).
\end{equation}
\end{lemma}

\begin{proof1}{\it of Theorem~{\rm \ref{prop:ctl1}.}} \ Note that $\de_{j}^{(2,m_j)}$ is the squared Euclidian norm of the centred Gaussian vector 
$\vec{\de}_j^{\,(m_j)}:=(\de_{j 1},\ldots,\de_{j m_j}).$
Therefore, $\de_{j}^{(2,m_j)}$ has the same distribution as $\sum_{k=1}^{m_j} \la_{jk}\,\ep_{jk}^2$, where $\la_{j 1},\ldots,\la_{j m_j}$ are the non-negative eigenvalues of the covariance matrix of $\vec{\de}_j^{\,(m_j)}$ and $\ep_{j 1},\ldots,\ep_{j m_j}$ are independent standard Gaussian random variables. Thus, using a version of the Lindeberg condition (see for instance \cite{CR} or Lemma~2 in \cite{IL}), it turns out that for proving the proposition it is enough to show that 
\begin{equation}
\label{prop:ctl1:eq4}
\lim_{j\rightarrow +\infty} \frac{\max_{1\le k \le m_j} \la_{jk}}{\sqrt{\var \left(\de_j^{(2,m_j)}\right)}}=0.
\end{equation}
To derive (\ref{prop:ctl1:eq4}) let us first prove that there is a positive constant $c_4$ (not depending on~$j$), such that  for all large enough $j$,
\begin{equation}
\label{prop:ctl1:eq5}
\var (\de_j^{(2,m_j)})\ge c_4 m_j\,.
\end{equation}
Using (\ref{prop:ctl1:eq1}), (\ref{varde}) and the change of variable $\eta=a_j\xi$, one gets
\begin{equation}
\label{vardel}
 \var(\de_j ^{(2,m_j)})=\sum_{k=1}^{m_j}\sum_{l=1}^{m_j}\cov(\de_{jk}^2,\de_{jl}^2)=2\sum_{k=1}^{m_j}\sum_{l=1}^{m_j}\cov^2(\de_{jk},\de_{jl})\ge 2\sum_{k=1}^{m_j} \var^2 (\de_{jk}) 
 \end{equation}
\begin{equation}
\label{prop:ctl1:eq6} =2m_j \Big (a_j\int_\R \frac{\big |\wh{\psi}(a_j\xi)\big |^2 h(\xi)}{|\xi^2-s_0^2|^{2\a}}\, d\xi\Big )^2=2m_j \Big (\int_\R \frac{\big |\wh{\psi}(\eta)\big |^2 h(a_j^{-1}\eta)}{|a_j^{-2}\eta^2-s_0^2|^{2\a}}\, d\eta\Big )^2\,.
\end{equation}
Moreover, it follows from (\ref{lem:var-de:eq2}) that 
\begin{equation}
\label{prop:ctl1:eq7}
\lim_{j\rightarrow +\infty} \int_\R \frac{\big |\wh{\psi}(\eta)\big |^2 h(a_j^{-1}\eta)}{|a_j^{-2}\eta^2-s_0^2|^{2\a}}\, d\eta=s_0^{-4\a} \int_\R \big |\wh{\psi}(\eta)\big |^2 \, d\eta >0.
\end{equation}
Then, (\ref{prop:ctl1:eq5}) results from (\ref{prop:ctl1:eq6}) and (\ref{prop:ctl1:eq7}).

Next, by Lemma~\ref{lem:sp-radius}  for all $j\in\N$ it holds
\begin{equation}
\label{prop:ctl1:eq9}
\max_{1\le k \le m_j} \la_{jk}\le \max_{1\le k \le m_j} \sum_{l=1}^{m_j}\big |\cov (\de_{jk},\de_{jl})\big |.
\end{equation}
Moreover, by (\ref{eq:bgam})  and (\ref{lem:cov-b:eq1}), for each fixed large enough $j$ and for every  $k\in\{1,\ldots, m_j\}$, one has 
 \[ \sum_{l=1}^{m_j}\big |\cov (\de_{jk},\de_{jl})\big |\le c_3\Big (1+a_j \sum_{l=1,\,l\ne k}^{m_j}|b_{jk}-b_{jl}|^{-1}\Big) \]
 \[\le c_3\Big (1+\frac{a_j}{\ga_j} \sum_{l=1,\,l\ne k}^{m_j}|k-l |^{-1}\Big)\le c_3\Big (1+\frac{2a_j}{\ga_j} \sum_{l=1}^{m_j}l^{-1}\Big) \]
\begin{equation}\label{prop:ctl1:eq10} \le c_3\Big (1+\frac{2a_j}{\ga_j} +\frac{2a_j}{\ga_j} \int_{1}^{m_j}y^{-1}\,dy\Big)\le c_3\bigg (1+\frac{2a_j\big (1+\log(m_j)\big)}{\ga_j}\bigg).
\end{equation}
Recall that the constant $c_3$ does not depend on $(j,k,l)$. Finally, putting together (\ref{prop:ctl1:eq2}), (\ref{prop:ctl1:eq5}), (\ref{prop:ctl1:eq9}), (\ref{prop:ctl1:eq10}),  and the fact that $\lim_{j\rightarrow +\infty} m_j=+\infty$, one gets (\ref{prop:ctl1:eq4}).
\end{proof1}

To obtain the exact asymptotic variance of $\de_j^{(2,m_j)}$ we specify asymptotic behaviours of the increments of the sequences $\{b_{jk}\}_{(j,k)\in \N\times \Z}$  in Assumption~\ref{Assumption_3}.
\newcounter{countD}
\renewcommand{\thecountD}{\arabic{countD}'}
\newtheorem{assumptionmD}[countD]{Assumption}

\setcounter{countD}{2}
\begin{assumptionmD}\label{Assumption_3'}
For all $j\in\N$ and for every $(k,l)\in\Z^2$ it holds 
\[
b_{jk}-b_{jl}= \ga_j (k-l),
\]
where $\{\ga_j\}_{j\in\N}$ is a sequence of positive real numbers such that 
\[\lim_{j\to +\infty}\frac{a_j}{\ga_j} = c \in (0,+\infty)\quad \mbox{and}\quad \lim_{j\to +\infty}m_j^2\left(\frac{\ga_j}{a_j}-\frac{1}{c}\right) = 0 .\]
\end{assumptionmD}
\begin{remark}
For example, Assumption~\ref{Assumption_3'} is satisfied for the sequence $\{\ga_j\}_{j\in\N}$ with  $\ga_j =a_j$ for all $j \ge j_0\in \mathbb{N}.$ 
\end{remark}

\begin{lemma}
\label{lem:sig-de2}
Let Assumption~{\rm\ref{Assumption_3'}} hold true and
\begin{equation}
\label{lem:sig-de2:eq0}
\lim_{j\rightarrow +\infty} m_j a_j^{-8}=0.
\end{equation}
Then, the sequence of positive real numbers $\big\{\var (\de_j^{(2,m_j)})/{m_j}\big\}_{j\in\N}$ converges to a finite and strictly positive limit when $j\rightarrow +\infty.$ More precisely, 
\begin{equation}
\label{lem:sig-de2:eq2}
\lim_{j\rightarrow +\infty} \frac{\var (\de_j^{(2,m_j)})}{m_j}=\V_1:=4c\pi s_0^{-8\a}  \int_{-c\pi} ^{c\pi} \Big |\sum_{n\in\Z} \big | \wh{\psi} (\eta+2nc\pi)\big|^2 \Big |^2\,d\eta .
\end{equation}
\end{lemma}

\begin{proof1}{\it of Lemma {\rm\ref{lem:sig-de2}.}} \ Using (\ref{eq:de}), (\ref{lem:var-de:eq1}), (\ref{vardel}),  Assumption~{\rm\ref{Assumption_3'}}, and the change of variable $\eta=a_j\xi$ one obtains
\[
\frac{\var (\de_j^{(2,m_j)})}{m_j}=\frac{2}{m_j}\sum_{k=1}^{m_j}\sum_{l=1}^{m_j}\I_{{\gamma_j}(k-l)/{a_j}}^2 (a_j^{-1}),\quad j\in\N,
\]
where $\I_\zeta^2(\cdot)$ is the squared function $\I_\zeta(\cdot)$ defined in (\ref{lem:var-de:eq1}). 

Let us denote by $F_j(\cdot)$ a bounded function defined on  $[-\pi a_j/\gamma_j,\pi a_j/\gamma_j]$ as
\[
F_j(\eta):= \sum_{n\in\Z} \big | \wh{\psi} (\eta+2n\pi a_j/\gamma_j)\big|^2.
\]
Let $\{\mu_j(k)\}_{k\in\Z}$ be the sequence of the Fourier coefficients of $F_j.$ These coefficients are real-valued since $\wh{\psi}(\cdot)$ is even. Using the fact that $\eta\mapsto e^{i {\gamma_j} k\eta/{a_j}} $ is, for each fixed $k\in\Z$, a $2\pi a_j/\gamma_j$-periodic function of $\eta$ and the dominated convergence theorem, one gets 
\begin{equation}
\label{lem:sig-de2:eq5}
\mu_j(k):=\int_{-\pi a_j/\gamma_j}^{\pi a_j/\gamma_j}e^{i {\gamma_j} k\eta/{a_j}} \Big (\sum_{n\in\Z} \big |\wh{\psi} (\eta+2n\pi a_j/\gamma_j)\big |^2\Big)\,d\eta=\int_\R e^{i {\gamma_j} k\eta/{a_j}} \big |\wh{\psi} (\eta)\big|^2\,d\eta.
\end{equation}
Now, let us  show that there is a finite constant $c_4$ such that, for all $j$ large enough, one has
 \[
m_j ^{-1/2}\bigg| \Big (\sum_{k=1}^{m_j}\sum_{l=1}^{m_j}\I_{{\gamma_j}(k-l)/{a_j}}^2 (a_j^{-1})\Big)^{1/2}-\Big ( \sum_{k=1}^{m_j}\sum_{l=1}^{m_j} s_0^{-8\a}\mu_j^2(k-l)\Big )^{1/2} \bigg |
\]
 \begin{equation}
\label{lem:sig-de2:eq6}\le 
c_4 \Big (m_j a_j^{-8}+a_j^{-4}\Big)^{1/2}.
\end{equation}
By the triangle inequality it holds  
\[
\bigg| \Big (\sum_{k=1}^{m_j}\sum_{l=1}^{m_j}\I_{{\gamma_j}(k-l)/{a_j}}^2 (a_j^{-1})\Big)^{1/2}-\Big ( \sum_{k=1}^{m_j}\sum_{l=1}^{m_j} s_0^{-8\a}\mu_j^2(k-l)\Big )^{1/2} \bigg |
\]
\begin{equation}
\label{lem:sig-de2:eq7}\le  \Big (\sum_{k=1}^{m_j}\sum_{l=1}^{m_j}
\big |\I_{{\gamma_j}(k-l)/{a_j}}(a_j^{-1})-s_0^{-4\a}\mu_j(k-l)\big|^2\Big)^{1/2}.
\end{equation}
Next, observe that it follows from (\ref{lem:var-de:eq2}), (\ref{lem:sig-de2:eq5}) and the inequalities $0<\a<1/2$ and $s_0>1$, that  for all $j$ large enough and for all $(k,l)\in\Z^2$ it holds
\[ \left|\I_{{\gamma_j}(k-l)/{a_j}}(a_j^{-1})-\frac{\mu_j(k-l)}{s_0^{4\a}}\right|^2\le\Bigg(\bigg | \I_{{\gamma_j}(k-l)/{a_j}} (a_j^{-1})-\frac{\int_{\R}e^{i{\gamma_j}(k-l)\eta/{a_j}}|\wh{\psi}(\eta)|^2\,d\eta}{s_0^{4\a}}\]
\[ -\frac{2\a}{s_0^{4\a-2}}\int_{\R}e^{i{\gamma_j}(k-l)\eta/{a_j}}\eta^2|\wh{\psi}(\eta)|^2\,d\eta \cdot a_j^{-2} \bigg| +\Big |\int_{\R}e^{i{\gamma_j}(k-l)\eta/{a_j}}\eta^2|\wh{\psi}(\eta)|^2\,d\eta\Big |a_j^{-2}\Bigg)^2\]
\[ \le 2 \bigg | \I_{{\gamma_j}(k-l)/{a_j}} (a_j^{-1})-s_0^{-4\a}\int_{\R}e^{i{\gamma_j}(k-l)\eta/{a_j}}|\wh{\psi}(\eta)|^2\,d\eta-2\a s_0^{-4\a-2}a_j^{-2}\]
\[ \times\int_{\R}e^{i{\gamma_j}(k-l)\eta/{a_j}}\eta^2|\wh{\psi}(\eta)|^2\,d\eta\bigg|^2+2\Big |\int_{\R}e^{i{\gamma_j}(k-l)\eta/{a_j}}\eta^2|\wh{\psi}(\eta)|^2\,d\eta\Big |^2 a_j^{-4}\]
\begin{equation}
\label{lem:sig-de2:eq7bis}
 \le 2 c_1^2 a_j^{-8}+2\Big |\int_{\R}e^{i{\gamma_j}(k-l)\eta/{a_j}}\eta^2|\wh{\psi}(\eta)|^2\,d\eta\Big |^2 a_j^{-4},
\end{equation} 
where $c_1$ is the constant from (\ref{lem:var-de:eq2}). 

By (\ref{lem:sig-de2:eq7}) and (\ref{lem:sig-de2:eq7bis}) to derive (\ref{lem:sig-de2:eq6}) it is sufficient to show that
\[
\sum_{k\in\Z}\Big |\int_{\R}e^{i {\gamma_j} k\eta/{a_j}}\eta^2|\wh{\psi}(\eta)|^2\,d\eta\Big |^2=
\sum_{k\in\Z}\Big |\int_{-\pi a_j/\gamma_j}^{\pi a_j/\gamma_j}e^{i {\gamma_j} k\eta/{a_j}}\sum_{n\in\Z} (\eta+2n\pi a_j/\gamma_j)^2 \]
\[\times\big | \wh{\psi} (\eta+2n\pi a_j/\gamma_j)\big|^2\,d\eta\Big |^2<+\infty.
\]
This inequality holds by Plancherel's identity as $\Big\{\int_{\R}e^{i {\gamma_j} k\eta/{a_j}}\eta^2|\wh{\psi}(\eta)|^2\,d\eta\Big\}_{k\in\Z}$ is the sequence of the Fourier coefficients of the bounded on $[-\pi a_j/\gamma_j,\pi a_j/\gamma_j]$ function 
$ \sum_{n\in\Z} (\eta+2n\pi a_j/\gamma_j)^2 \big | \wh{\psi} (\eta+2n\pi a_j/\gamma_j)\big|^2.$

Next, let us define $F_0(\cdot)$  as
\begin{equation}
\label{lem:sig-de2:eq4}
F_0(\eta):=\sum_{n\in\Z} \big |\wh{\psi} (\eta+2nc\pi)\big |^2, \quad \eta\in [-c\pi, c\pi],
\end{equation}
where $c$ is the same positive constant as in Assumption~{\rm\ref{Assumption_3'}}. $F_0(\cdot)$ is a bounded function on $[-c\pi, c\pi].$  

Let us now show that 
\begin{equation}
\label{lem:sig-de2:eq9} 
\lim_{j\rightarrow +\infty} \frac{1}{m_j}\sum_{k=1}^{m_j}\sum_{l=1}^{m_j} \mu_j^2(k-l)=2c\pi\int_{-c\pi} ^{c\pi} |F_0(\eta)|^2\,d\eta.
\end{equation} 
Note that \[ \frac{1}{m_j}\sum_{k=1}^{m_j}\sum_{l=1}^{m_j} \mu_j^2(k-l)=\frac{1}{m_j}\sum_{k=1}^{m_j}\sum_{q=k-m_j}^{k-1}\mu_j^2(q)\]
and for the sequence $\{\mu_0(k)\}_{k\in\Z}$ of the Fourier coefficients of $F_0$ it holds
\begin{equation}
\label{diff} \frac{1}{m_j}\left|\sum_{k=1}^{m_j}\sum_{q=k-m_j}^{k-1}\mu_j^2(q)-\sum_{k=1}^{m_j}\sum_{q=k-m_j}^{k-1}\mu_0^2(q)\right|\le \frac{C}{m_j}\sum_{k=1}^{m_j}\sum_{q=k-m_j}^{k-1}\left|\mu_j(q)-\mu_0(q)\right|\end{equation} 
as $\mu_j(q)$ and $\mu_0(q)$ are bounded by $\int_\R \big |\wh{\psi} (\eta)\big|^2\,d\eta.$

Using the expressions for Fourier coefficients and Assumption~\ref{Assumption_2},  we get that for $k=1,...,m_j$ 
\[\sum_{q=k-m_j}^{k-1}\left|\mu_j(q)-\mu_0(q)\right|\le \sum_{q=k-m_j}^{k-1}\int_{-A}^A \left|e^{i \frac{{\gamma_j} k\eta}{{a_j}}}- e^{i \frac{ k\eta}{c}} \right| \big |\wh{\psi} (\eta)\big|^2\,d \eta\]
\[\le C' \sum_{q=k-m_j}^{k-1}\int_{-A}^A \left|\sin \left(\frac{k\eta}{2}\left(\frac{{\gamma_j} }{a_j}-\frac{1}{c}\right)\right)\right|\,d \eta\]
\[ \le C' \sum_{q=-m_j}^{m_j}\int_{-A}^A \left|\sin \left(\frac{k\eta}{2}\left(\frac{{\gamma_j} }{a_j}-\frac{1}{c}\right)\right)\right|\,d \eta. \]
Hence, it follows from the inequality $|\sin(x)|\le |x|$ and Assumption~\ref{Assumption_3'} that
\begin{equation}
\label{boundsum}
\sum_{q=k-m_j}^{k-1}\left|\mu_j(q)-\mu_0(q)\right|\le C'' m_j^2\left|
\frac{{\gamma_j} }{a_j}-\frac{1}{c}\right|\to 0, \ j\to +\infty.
\end{equation} 
Thus, by (\ref{diff}), (\ref{boundsum}) and the Ces\`{a}ro mean convergence theorem one gets
\begin{equation}
\label{sumdif} \frac{1}{m_j}\left|\sum_{k=1}^{m_j}\sum_{q=k-m_j}^{k-1}\mu_j^2(q)-\sum_{k=1}^{m_j}\sum_{q=k-m_j}^{k-1}\mu_0^2(q)\right|\to 0, \ j\to 0.
\end{equation} 
Now, by Plancherel's identity 
\[\frac{1}{m_j}\sum_{k=1}^{m_j}\sum_{q=k-m_j}^{k-1}\mu_0^2(q)=\sum_{q=-\infty}^{+\infty} \mu_0^2(q)-\frac{1}{m_j}\sum_{k=1}^{m_j}\sum_{q=k}^{+\infty} \mu_0^2(q) -\frac{1}{m_j}\sum_{k=1}^{m_j}\,\sum_{q=-\infty}^{k-m_j-1}\mu_0^2(q)\]
\begin{equation}
\label{lem:sig-de2:eq10} =2c\pi\int_{-c\pi}^{c\pi} \left|F_0(\eta)\right|^2\,d\eta
 -\frac{1}{m_j}\sum_{k=1}^{m_j}\sum_{q=k}^{+\infty} \mu_0^2(q)-\frac{1}{m_j}\sum_{k'=1}^{m_j}\,\sum_{q=-\infty}^{-k'}\mu_0^2(q).
\end{equation}
Next, observe that the sequence $\big \{\sum_{q=k}^{+\infty} \mu_0^2(q)\big\}_{k\in\N}$ converges to zero. Consequently by the Ces\`{a}ro mean convergence theorem one gets 
\begin{equation}
\label{lem:sig-de2:eq11}
\lim_{j\rightarrow +\infty} \frac{1}{m_j}\sum_{k=1}^{m_j}\sum_{q=k}^{+\infty} \mu_0^2(q)=0.
\end{equation}
Using the same arguments, one obtains that
\begin{equation}
\label{lem:sig-de2:eq12}
\lim_{j\rightarrow +\infty}\frac{1}{m_j}\sum_{k'=1}^{m_j}\,\sum_{q=-\infty}^{-k'}\mu_0^2(q)=0.
\end{equation}
Putting together (\ref{sumdif}), (\ref{lem:sig-de2:eq10}), (\ref{lem:sig-de2:eq11}) and (\ref{lem:sig-de2:eq12}) it follows that 
(\ref{lem:sig-de2:eq9}) holds true.

Finally, combining (\ref{lem:sig-de2:eq9}) with (\ref{lem:sig-de2:eq0}), (\ref{lem:sig-de2:eq6}) and (\ref{lem:sig-de2:eq4}) one obtains (\ref{lem:sig-de2:eq2}).
\end{proof1}

\section{Asymptotic normality of two auxiliary statistics}\label{sec_4}

This section proves asymptotic normality of two auxiliary statistics of the semiparametric model defined by Assumption~\ref{Assumption_1}. They are two functions of the parameters $s_0$ and $\alpha.$ The results will be used in the following sections to derive and investigate simultaneous estimators of  $s_0$ and $\alpha.$

Let us set 
\begin{equation}
\label{theo:ctl1:eq2}
\ovl{\de}_j ^{(2,m_j)}:=\frac{\de_j^{(2,m_j)}}{m_j}=\frac{1}{m_j}\sum_{k=1}^{m_j} \de_{jk}^2, \quad j\in\N,
\end{equation}
where $\de_{jk}$ is given in Definition~\ref{def2}.

The following theorem introduces the first statistics and derives its asymptotic normality. 

\begin{theorem}
\label{theo:ctl1}
Let the array $\{b_{jk}\}_{(j,k)\in\N\times\Z}$ satisfy Assumption~{\rm\ref{Assumption_3'}} and 
\begin{equation}
\label{theo:ctl1:eq1}
\lim_{j\rightarrow +\infty} m_j a_j ^{-4}=0.
\end{equation} 

Then, when $j$ goes to $+\infty$, the random variables
\begin{equation}
\label{theo:ctl1:eq3}
\ovl{Y}_j:=\sqrt{m_j} \left(\ovl{\de}_j ^{(2,m_j)}-s_0^{-4\a}\int_{\R}|\wh{\psi}(\eta)|^2\,d\eta\right)
\end{equation} 
converge in distribution to a centred Gaussian random variable $\ovl{Y}$ with the variance $\var (\ovl{Y})=\V_1$ given by {\rm (\ref{lem:sig-de2:eq2}).}
\end{theorem}
\begin{remark}
If the array $\{b_{jk}\}_{(j,k)\in\N\times\Z}$ satisfies Assumption~\ref{Assumption_3'}, then the condition (\ref{prop:ctl1:eq2}) of Theorem~\ref{prop:ctl1} holds true for any $\{m_j\}_{j\in\N}$.
\end{remark}

\begin{proof1}{\it of Theorem~{\rm \ref{theo:ctl1}.}}\ 
By Theorem~\ref{prop:ctl1}, when $j$ goes to $+\infty$, the random variables $\sqrt{\V_1}\, Y_j$ converge in distribution to a centred Gaussian random variable $\ovl{Y}$ whose variance equals $\V_1.$ Moreover, by (\ref{prop:ctl1:eq3}) and (\ref{theo:ctl1:eq2}) the random variable $\sqrt{\V_1}\,Y_j$ equals
\[
\sqrt{\V_1}\,Y_j=\sqrt{\V_1\times\frac{m_j^{}}{\var\left(\de_j^{(2,m_j)}\right)}}\,\sqrt{m_j}\, \Big (\ovl{\de}_j ^{(2,m_j)}-\E\left(\ovl{\de}_j ^{(2,m_j)}\right)\Big),
\]
and, by Lemma~\ref{lem:sig-de2}, it holds 
\[
\lim_{j\rightarrow +\infty} \sqrt{\V_1\times\frac{m_j^{}}{\var\left(\de_j^{(2,m_j)}\right)}}=1.
\]
Thus, when $j$ goes to $+\infty,$ the random variables $\sqrt{m_j} \Big (\ovl{\de}_j ^{(2,m_j)}-\E\left(\ovl{\de}_j ^{(2,m_j)}\right)\Big)$ converge in distribution to $\ovl{Y}.$   To show that the sequence $\big \{\ovl{Y}_j \big\}_{j\in\N}$ shares the same property, it is enough to prove that 
\begin{equation}
\label{theo:ctl1:eq6}
\lim_{j\rightarrow +\infty} \sqrt{m_j} \left(\E\left(\ovl{\de}_j ^{(2,m_j)}\right)-s_0^{-4\a}\int_{\R}|\wh{\psi}(\eta)|^2\,d\eta\right)=0.
\end{equation}
It follows from from (\ref{varde}), (\ref{lem:var-de:eq1}) and (\ref{theo:ctl1:eq2})  that 
$\E\left(\ovl{\de}_j ^{(2,m_j)}\right)=\I_0 (a_j^{-1})$. Thus, using Lemma~\ref{lem:var-de} and (\ref{theo:ctl1:eq1}) one obtains (\ref{theo:ctl1:eq6}).
\end{proof1}

Let $\{M_j\}_{j\in\N}$ be a sequence of positive integers defined as
\begin{equation}
\label{theo:ctl2:eq4}
M_j:=\left[\frac{m_j}{(a^{-2}_{j+1}-a^{-2}_{j+2})^{2}}\right]\,,
\end{equation}
where $[\cdot]$ denotes the integer part function. 

\begin{remark}\label{Mj_increase}
$\{M_j\}_{j\in\N}$ is an increasing sequence as by (\ref{theo:ctl2:eq4}) one gets
 \[
 M_j=\left[\frac{m_j(a_{j+1}a_{j+2})^{4}}{(a^{2}_{j+2}-a^{2}_{j+1})^{2}} \right]\ge\left[m_ja_{j+1}^{4}\right]\to +\infty\,.
 \]
\end{remark}


\newcounter{countB}
\renewcommand{\thecountB}{\arabic{countB}*}
\newtheorem{assumptionmB}[countB]{Assumption}

\setcounter{countB}{2}
\begin{assumptionmB}\label{Assumption_3*}
For all $j\in\N$ and for every $(k,l)\in\Z^2$ it holds 
\[
b_{jk}-b_{jl}= \ga_j (k-l),
\]
where $\{\ga_j\}_{j\in\N}$ is a sequence of positive real numbers such that 
\[\lim_{j\to +\infty}\frac{a_j}{\ga_j} = c \in (0,+\infty)\quad \mbox{and}\quad \lim_{j\to +\infty}m_j^2a_j^8\left(\frac{\ga_j}{a_j}-\frac{1}{c}\right) = 0 .\]
\end{assumptionmB}
\begin{remark}
For example, Assumption~\ref{Assumption_3*} is satisfied if for all $j \ge j_0\in \mathbb{N}$ it holds $\ga_j =a_j.$  
\end{remark}

 Now we introduce the second auxiliary statistics \[
\De\ovl{\de}_{j+1}^{(2,M_j)}:=\frac{\ovl{\de}_{j+1}^{(2,M_j)}-\ovl{\de}_{j+2}^{(2,M_j)}}{a_{j+1}^{-2}-a_{j+2}^{-2}}
\] via increments of $\ovl{\de}_j ^{(2,M_j)}$ and prove its asymptotic normality. 
\begin{theorem}
\label{theo:ctl2}
Assume that the following conditions hold: 
\begin{enumerate}
\item There exists $B\in (0,A)$ such that $\wh{\psi}$ vanishes on the interval $[-B,B]$, that is 
\begin{equation}
\label{theo:ctl2:eq1}
\sp\,\wh{\psi}\subseteq \big\{\xi\in\R: B\le |\xi|\le A\big\}.
\end{equation}
\item Assumption~{\rm \ref{Assumption_3*}} holds true and for some $j_0\in\N$ the sequence 
$\{a_j\}_{j\in\N}$ satisfies
\begin{equation}
\label{theo:ctl2:eq3}
\frac{ a_{j+1}}{a_j}\ge \frac{A}{B}>1\,,\quad \mbox{for all}\quad j\ge j_0.
\end{equation}

\item The sequence $\{m_j\}_{j\in\N}$ satisfies {\rm (\ref{theo:ctl1:eq1})}.
\end{enumerate}
 Then, when $j$ goes to $+\infty$, the random variables
\begin{equation}
\label{theo:ctl2:eq5} 
\ovl{Z}_j:=\sqrt{m_j}\bigg (\De\ovl{\de}_{j+1}^{(2,M_j)}
-2\a s_0^{-4\a-2}\int_{\R}\eta^2|\wh{\psi}(\eta)|^2\,d\eta\bigg)
\end{equation}
converge in distribution to a centred Gaussian random variable $\ovl{Z}$ with the variance $\var (\ovl{Z})=2\V_1$.
\end{theorem}
\begin{remark}\label{rem5}
Notice that (\ref{theo:ctl2:eq1}) and (\ref{theo:ctl2:eq3}) imply that $\sp\,\wh{\psi}(a_j \cdot)\bigcap\,\sp\,\wh{\psi}(a_{j+1} \cdot)$ is a Lebesgue negligible set for all sufficiently large $j\in\N.$
\end{remark}

\begin{proof1}{\it of Theorem~{\rm\ref{theo:ctl2}.}}\ First notice that it follows from (\ref{eq:de}) and Remark~\ref{rem5} that $\cov (\de_{(j+1)k},\de_{(j+2)l})=0$ for all  $(k,l)\in\{1,\ldots, M_j\}^2$ and sufficiently large $j\in\N,$ which means that the centred Gaussian vectors $\vec{\de}_{j+1}^{\, (M_j)}:=(\de_{(j+1)1},\ldots,\de_{(j+1)M_j})$ and $\vec{\de}_{j+2}^{\, (M_j)}:=(\de_{(j+2)1},\ldots,\de_{(j+2)M_j})$ are independent. Therefore, the two random variables
\[ 
 \de_{j+1}^{(2,M_j)}:=\sum_{k=1}^{M_j} \de_{(j+1)k}^2
\quad \mbox{and} \quad
 \de_{j+2}^{(2,M_j)}:=\sum_{k=1}^{M_j} \de_{(j+2)k}^2
\]
are independent. 

By Remark~\ref{Mj_increase} $\{M_j\}_{j\in\N}$ is an increasing sequence. Hence, by Assumption~\ref{Assumption_3*} condition (\ref{prop:ctl1:eq2}) is satisfied if $m_j$ is replaced by $M_{j-1}$ or by $M_{j-2}$. Therefore, by Theorem~\ref{prop:ctl1}, when $j$ goes to $+\infty,$ the random variables 
\[
Z_{1,j}:=\frac{\de_{j+1}^{(2,M_j)}-\E (\de_{j+1}^{(2,M_j)})}{\sqrt{\var \left(\de_{j+1}^{(2,M_j)}\right)}}
\]
converge in distribution to a standard Gaussian random variable, and that the random variables 
\[
Z_{2,j}:=\frac{\de_{j+2}^{(2,M_j)}-\E (\de_{j+2}^{(2,M_j)})}{\sqrt{\var \left(\de_{j+2}^{(2,M_j)}\right)}}
\]
share the same property. 

Next, using (\ref{theo:ctl1:eq1}), (\ref{theo:ctl2:eq4})  and (\ref{theo:ctl2:eq3}), one gets that 
\[\lim_{j\rightarrow +\infty} \frac{M_j}{a^{8}_{j+1}}=\lim_{j\rightarrow +\infty}\left (\frac{m_j}{a_{j+1}^{4}}\cdot\frac{\left({a_{j+2}}/{a_{j+1}}\right)^{4}}{\left(\left({a_{j+2}}/{a_{j+1}}\right)^{2}-1\right)^{2}} \right)=0\]
as the function $\frac{x^4}{(x^2-1)^2}$ is bounded from above for $x\in [A/B,+\infty).$ The same is also true for ${M_j}/{a^{8}_{j+2}}.$

 Therefore, by Lemma~\ref{lem:sig-de2} 
\[
\lim_{j\rightarrow +\infty} \frac{\sqrt{\var \left(\de_{j+1}^{(2,M_j)}\right)}}{\sqrt{M_j}}=\sqrt{\V_1} \quad\mbox{and}\quad 
\lim_{j\rightarrow +\infty} \frac{\sqrt{\var \left(\de_{j+2}^{(2,M_j)}\right)}}{\sqrt{M_j}}=\sqrt{\V_1}.
\]
Thus, when $j$ goes to $+\infty$, the sequence 
\[
Z_{1,j}':=\frac{\sqrt{\var \left(\de_{j+1}^{(2,M_j)}\right)}}{\sqrt{M_j}}\, Z_{1,j}=\frac{\de_{j+1}^{(2,M_j)}-\E (\de_{j+1}^{(2,M_j)})}{\sqrt{M_j}}
\]
converges in distribution to a centred Gaussian random variable with variance $\V_1$, and the sequence 
\[
Z_{2,j}':=\frac{\sqrt{\var \left(\de_{j+2}^{(2,M_j)}\right)}}{\sqrt{M_j}}\, Z_{2,j}=\frac{\de_{j+2}^{(2,M_j)}-\E (\de_{j+2}^{(2,M_j)})}{\sqrt{M_j}}
\]
shares the same property. Therefore, using the fact that for sufficiently large $j$ these two sequences are independent and the equalities $\E (\de_{j+1}^{(2,M_j)})=M_j \I_0 (a^{-1}_{j+1})$ and $\E (\de_{j+2}^{(2,M_j)})=M_j \I_0 (a^{-1}_{j+2})$, one gets that the random variables
\begin{eqnarray*}
Z_{1,j}'-Z_{2,j}'&=&\frac{\de_{j+1}^{(2,M_j)}-\de_{j+2}^{(2,M_j)}}{\sqrt{M_j}}-\sqrt{M_j}\big (\I_0 (a^{-1}_{j+1})-\I_0 (a^{-1}_{j+2})\big)\\
&=& \sqrt{M_j} \Big (\ovl{\de}_{j+1}^{(2,M_j)}-\ovl{\de}_{j+2}^{(2,M_j)}-\big (\I_0 (a^{-1}_{j+1})-\I_0 (a^{-1}_{j+2})\big)\Big)
\end{eqnarray*}
converge in distribution to a centred Gaussian random variable with the variance $2\V_1,$ when $j\to +\infty.$ 

By (\ref{theo:ctl2:eq4}) the sequence of
\begin{eqnarray}
\ovl{Z}_j '&:=& \frac{\sqrt{m_j}\big (a^{-2}_{j+1}-a^{-2}_{j+2})^{-1}}{\sqrt{M_j}}\big (Z_{1,j}'-Z_{2,j}'\big)\nonumber\\
&=&\sqrt{m_j}\Bigg (\frac{\ovl{\de}_{j+1}^{(2,M_j)}-\ovl{\de}_{j+2}^{(2,M_j)}}{a^{-2}_{j+1}-a^{-2}_{j+2}}
-\frac{\I_0 (a^{-1}_{j+1})-\I_0 (a^{-1}_{j+2})}{a^{-2}_{j+1}-a^{-2}_{j+2}}\Bigg)\nonumber
\end{eqnarray}
shares the same property. 

Thus, it turns out that for deriving the theorem it is enough to show that
\begin{equation}
\label{theo:ctl2:eq10}
\lim_{j\rightarrow +\infty} \sqrt{m_j}\left(\frac{\I_0 (a^{-1}_{j+1})-\I_0 (a^{-1}_{j+2})}{a^{-2}_{j+1}-a^{-2}_{j+2}}
-2\a s_0^{-4\a-2}\int_{\R}\eta^2|\wh{\psi}(\eta)|^2\,d\eta\right)=0.
\end{equation}
Using Lemma~\ref{lem:var-de} one gets that 
\[\bigg  | \I_0 (a^{-1}_{j+1})-\I_0 (a^{-1}_{j+2})-\Big(2\a s_0^{-4\a-2}\int_{\R}\eta^2|\wh{\psi}(\eta)|^2\,d\eta\Big)(a^{-2}_{j+1}-a^{-2}_{j+2})\bigg |\]
\[\le \bigg | \I_0 (a^{-1}_{j+1})-s_0^{-4\a}\int_{\R}|\wh{\psi}(\eta)|^2\,d\eta-\Big(2\a s_0^{-4\a-2}\int_{\R}\eta^2|\wh{\psi}(\eta)|^2\,d\eta\Big)a^{-2}_{j+1}\bigg|\]
\[+\bigg | \I_0 (a^{-1}_{j+2})-s_0^{-4\a}\int_{\R}|\wh{\psi}(\eta)|^2\,d\eta-\Big(2\a s_0^{-4\a-2}\int_{\R}\eta^2|\wh{\psi}(\eta)|^2\,d\eta\Big)a^{-2}_{j+2}\bigg|\]
\[
\le c_1 \big(a^{-4}_{j+1}+a^{-4}_{j+2}\big),\]
where $c_1$ is the constant in (\ref{lem:var-de:eq2}). Thus, 
\begin{equation}
\label{theo:ctl2:eq11} 
\sqrt{m_j}\bigg |\frac{\I_0 (a^{-1}_{j+1})-\I_0 (a^{-1}_{j+2})}{a^{-2}_{j+1}-a^{-2}_{j+2}}
-2\a s_0^{-4\a-2}\int_{\R}\eta^2|\wh{\psi}(\eta)|^2\,d\eta\bigg|\le \frac{c_1 \sqrt{m_j}(a^{-4}_{j+1}+a^{-4}_{j+2})}{a^{-2}_{j+1}-a^{-2}_{j+2}}.
\end{equation}
Finally, combining (\ref{theo:ctl1:eq1}), (\ref{theo:ctl2:eq3}) and  (\ref{theo:ctl2:eq11})  one gets
\[\frac{\sqrt{m_j}(a^{-4}_{j+1}+a^{-4}_{j+2})}{a^{-2}_{j+1}-a^{-2}_{j+2}}=\frac{\sqrt{m_j}}{a_{j+1}^2}\cdot\frac{1+\left({a_{j+1}}/{a_{j+2}}\right)^4}{1-\left({a_{j+1}}/{a_{j+2}}\right)^2}\to 0,\ j\to +\infty,\]
which confirms (\ref{theo:ctl2:eq10}) and finishes the proof.
\end{proof1}

\begin{remark}
For example, the sequence $\{a_j\}_{j\in \N}$ with $a_j=a^j,$ $j\in \N,$ and $a\ge A/B$ satisfies the assumptions of Theorem~\ref{theo:ctl2}.
\end{remark}

Note that under the conditions of Theorem~\ref{theo:ctl2}, for sufficiently large $j\in\N,$ the random variable $\ovl{Y}_j$ defined in (\ref{theo:ctl1:eq3}) is independent of $\ovl{Z}_j$ defined by (\ref{theo:ctl2:eq5}). It is easy to see as the centred Gaussian random vectors $\vec{\de}_{j}^{\, (m_j)}:=(\de_{j 1},\ldots,\de_{j m_j})$, $\vec{\de}_{j+1}^{\, (M_j)}:=(\de_{(j+1)1},\ldots,\de_{(j+1)M_j})$ and $\vec{\de}_{j+2}^{\, (M_j)}:=(\de_{(j+2)1},\ldots,\de_{(j+2)M_j})$ are independent.
Therefore, the following result follows from Theorems~\ref{theo:ctl1} and \ref{theo:ctl2}. 

\begin{corollary}
\label{cor:ctl-vec1}
When $j$ goes to $+\infty$, the  random vectors $(\ovl{Y}_j, \ovl{Z}_j)$ converge in distribution to the random vector $(\ovl{Y},\ovl{Z})$ with the bivariate centred Gaussian distribution  
$
\cal N \left (\left (\begin{array}{c}
 0\\
0 \end{array}\right),
\left (\begin{array}{cc}
\V_1 & 0\\
0 & 2\V_1
\end{array}\right)
\right).
$
\end{corollary}

\section{Asymptotic normality of adjusted estimators}\label{sec_5}

In this section the axillary statistics $\ovl{\de}_j ^{(2,m_j)}$ and $\De\ovl{\de}_{j+1}^{(2,M_j)}$ are used for deriving adjusted statistics to estimate the parameters of interest. The central limit theorem is proved for the proposed adjusted statistics.

By (\ref{theo:ctl1:eq3}), (\ref{theo:ctl2:eq5}) and Corollary \ref{cor:ctl-vec1}, under the assumptions of Theorem~\ref{theo:ctl2} one has
\begin{equation}
\label{theo:ctl3:eq_clt1et2} 
 \sqrt{m_j} \begin{pmatrix} \ovl{\de}_j ^{(2,m_j)}-s_0^{-4\a}\int_{\R}|\wh{\psi}(\eta)|^2\,d\eta \\~\\
                                \De\ovl{\de}_{j+1}^{(2,M_j)} -2\a s_0^{-4\a-2}\int_{\R}\eta^2|\wh{\psi}(\eta)|^2\,d\eta
                                \end{pmatrix} 
\mathop{\xrightarrow{\, d\, }} \mathcal{N}\left( 0, \begin{pmatrix} \V_1 & 0 \\ 0 & 2\V_1 \end{pmatrix}\right),
\end{equation}
when $\ {j \to +\infty}.$

This two-dimensional central limit theorem gives the fluctuation rate for the corresponding law of large number proven in \cite{AAFO}
\begin{equation}
\label{theo:ctl3:eq_lfgn} 
\left( \frac{ \ovl{\de}_j ^{(2,m_j)} }{ \int_{\R}|\wh{\psi}(\eta)|^2\,d\eta } ,
                \frac{ \De\ovl{\de}_{j+1}^{(2,M_j)} }{ 2\int_{\R}\eta^2|\wh{\psi}(\eta)|^2\,d\eta }
                \right) 
\mathop{\xrightarrow{\,a.s.\,}} \Phi(s_0,\a):=
\left(s_0^{-4\a} ,
                \a s_0^{-4\a-2} \right),
\end{equation}
when $\ {j \to +\infty}.$

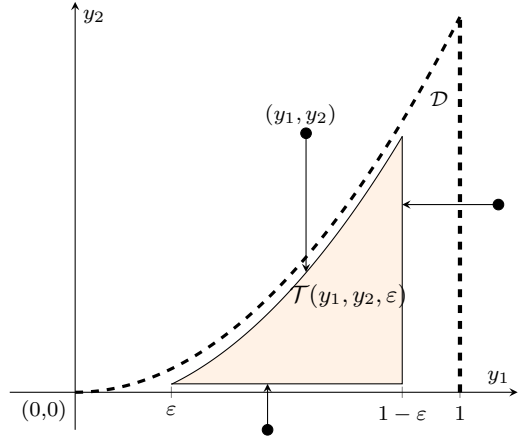
\begin{figure}[!htb]
\centering
\begin{tikzpicture}
\begin{axis}[width=0.7\textwidth,
       height=0.6\textwidth,xmax=1.15,xmin=-0.2, ymax=0.52, ymin=-0.05,
xtick={0,0.25,0.85,1},xticklabels={$0$,$\varepsilon$,$1-\varepsilon$,$1$},ytick={0,0}, x axis line style={draw=none},
   xlabel=$y_1$,axis lines=middle, ylabel=$y_2$
]

\addplot [name path=f,dashed,very thick,
    domain=0:1, 
    samples=100, 
    ]
{(x^2)/2};

\addplot [name path=g, domain=0.25:0.85, 
    samples=100, 
    ]
{(x^2)/2-0.020625};

\addplot [const plot] coordinates {(0.25,0.01125) (0.85,0.01125) (0.85,0.34125)};

\addplot [const plot] coordinates {(0,0) (1,0)};

 \path[name path=axis] (axis cs:0.25,0.01125) -- (axis cs:0.85,0.01125);
\addplot [
        thick,
        color=blue,
        fill=orange, 
        fill opacity=0.1
    ]
fill between[
        of= g and axis,
        soft clip={domain=0:1},
    ];
\addplot [const plot, dashed,ultra thick,
] coordinates {(1,0) (1,0.5)};

\draw[->,>=stealth] (axis cs:1,0) -- (axis cs:1.15,0);
\addplot[mark=*] coordinates {(0.5,-0.05)};
\draw[->,>=stealth] (axis cs:0.5,-0.05) -- (axis cs:0.5,0.01125);
\addplot[mark=*] coordinates {(1.1, 0.25)};
\draw[->,>=stealth] (axis cs:1.1, 0.25) -- (axis cs:0.85,0.25);
\addplot[mark=*] coordinates {(0.6, 0.345)};
\draw[->,>=stealth] (axis cs:0.6, 0.345) -- (axis cs:0.6, 0.160);
\draw[-] (axis cs:-0.17, 0) -- (axis cs:0,0);
\node at (axis cs:-0.0,-0.0) [anchor=north east] {(0,0)};

\node at (axis cs:0.99, 0.37) [anchor=south east] {$\mathcal{D}$} ;
\node at (axis cs:0.7, 0.345) [anchor=south east] {$\left(y_1, y_2\right)$} ;
\node at (axis cs:0.71, 0.155) [anchor=north] {\small{$\mathcal{T}({y_1, y_2, \varepsilon})$}} ;
\end{axis}
\end{tikzpicture}
\caption{Plot of ($y_1$, $y_2$) and the corresponding truncated values}\label{fig22}
\end{figure}

Let us consider the function $g:[-1,+\infty) \to [-1/e,+\infty)$ defined as $g(t)=te^t.$ This is an increasing continuous one-to-one function. Its inverse function is $\text{LambertW}$ that is continuous, defined on  $[-1/e,+\infty)$ with values in $[-1,+\infty)$ and satisfies 
\[ 
\text{LambertW}(y)~ e^{\text{LambertW}(y)} = y \quad\quad \text{ i.e. } \quad\quad e^{\text{LambertW}(y)} = \frac{y}{\text{LambertW}(y)}. 
\]
As stated in \cite{AAFO}, the vector-valued function
$\Phi:  (1,+\infty)\times(0,1/2) \to \mathcal{D}$ defined in (\ref{theo:ctl3:eq_lfgn}) is  a continuous one-to-one function taking values in \[\mathcal{D}=\left\{(y_1,y_2)\in\R^2~:~0<y_1<1~\text{ and } 0<y_2<\frac{y_1^2}{2} \right\}.\]
Its inverse function $\Phi^{-1}: \mathcal{D} \to (1,+\infty)\times(0,1/2)$ is continuous and given by 
\begin{align*} \Phi^{-1}(y_1,y_2)=& \left( \exp\left( \frac12 \text{LambertW}\left( -\frac{y_1\ln(y_1)}{2y_2} \right) \right),\right. \\
&\left.\frac{y_2}{y_1} \exp\left( \text{LambertW}\left( -\frac{y_1\ln(y_1)}{2y_2} \right) \right) \right).
\end{align*}
Let us define the following continuous vector-valued truncating function $\mathcal{T}$ defined for $\eps\in(0,1),$ $(y_1,y_2)\in\R^2,$ and taking values in $\mathcal{D}$
\[\mathcal{T}(y_1,y_2,\eps) = \Big( \mathcal{T}_1(y_1,\eps) ~,~ \mathcal{T}_2(y_1,y_2,\eps) \Big) \in \mathcal{D},\]
where
\begin{align*}  \mathcal{T}_1(y_1,\eps) &:= \max( \eps , \min( y_1, 1-\eps) ) = \begin{cases} 
                                                                   \eps,   & \text{ if } y_1 \le \eps, \\
                                                                   y_1,    & \text{ if } \eps \le y_1 \le 1-\eps, \\
                                                                   1-\eps, & \text{ if } y_1 > 1-\eps,
                                                                 \end{cases}  \\
  \mathcal{T}_2(y_1,y_2,\eps) &:= \max\left( {\eps^2}/{4} , \min\left( y_2, \frac{ \big(\mathcal{T}_1(y_1,\eps)\big)^2 }{2}-{\eps^2}/{4} \right) \right)\\
                          &\ \ = \begin{cases} 
                                   {\eps^2}/{4},   & \text{ if } y_2 \le {\eps^2}/{4}, \\
                                   y_2,    & \text{ if } {\eps^2}/{4} \le y_2 \le \frac{ \big(\mathcal{T}_1(y_1,\eps)\big)^2 }{2}-{\eps^2}/{4}, \\
                                   \frac{ \big(\mathcal{T}_1(y_1,\eps)\big)^2 }{2}-{\eps^2}/{4}, 
                                          & \text{ if } y_2 > \frac{ \big(\mathcal{T}_1(y_1,\eps)\big)^2 }{2}-{\eps^2}/{4}.
                                 \end{cases}  
\end{align*}
For values outside the feasible region $\mathcal{D},$ some typical mappings  by the truncating function $\mathcal{T}$ are sketched in Figure~\ref{fig22}.

Note that for each $(y_1,y_2)\in\mathcal{D}$   there is a small enough $\eps>0$ such that $\mathcal{T}(y_1,y_2,\eps)=(y_1,y_2)$ because $\mathcal{D}$ is an open set. Assumption~\ref{Assumption_1} on the parameters ensures that $(s_0,\a) \in (1,+\infty)\times(0,1/2)$ and therefore $\Phi(s_0,\a) \in \mathcal{D}.$ 

\begin{definition}
The adjusted statistic for the parameter $(s_0,\a)$ is  
\[ 
\widehat{(s_0,\a)}_j := \Phi^{-1} \left( 
                                  \mathcal{T}\Bigg( \frac{ \ovl{\de}_j ^{(2,m_j)} }{ \int_{\R}|\wh{\psi}(\eta)|^2\,d\eta } ~,~
                                                    \frac{ \De\ovl{\de}_{j+1}^{(2,M_j)} }{ 2\int_{\R}\eta^2|\wh{\psi}(\eta)|^2\,d\eta } ~,~
                                                    \frac{1}{m_j} \Bigg) \right).
\]
\end{definition}

Note that for some observations the values  $\left( \frac{ \ovl{\de}_j ^{(2,m_j)} }{ \int_{\R}|\wh{\psi}(\eta)|^2\,d\eta } ,
                \frac{ \De\ovl{\de}_{j+1}^{(2,M_j)} }{ 2\int_{\R}\eta^2|\wh{\psi}(\eta)|^2\,d\eta }
                \right)$ may not be in the feasible region $\mathcal{D}.$ Therefore, the truncation $\mathcal{T}$ was needed to guarantee that $\Phi^{-1}$ acts only on values from $\mathcal{D}.$ 

\begin{remark}
As for sufficiently large $j$ the vector $\left( \frac{ \ovl{\de}_j ^{(2,m_j)} }{ \int_{\R}|\wh{\psi}(\eta)|^2\,d\eta } ,
                \frac{ \De\ovl{\de}_{j+1}^{(2,M_j)} }{ 2\int_{\R}\eta^2|\wh{\psi}(\eta)|^2\,d\eta }
                \right)$ falls in $\mathcal{D},$ then $\widehat{(s_0,\a)}_j$ and  the corresponding adjusted statistic in \cite{AAFO} coincide almost surely. At the same time the new statistic requires only the simple truncation $\mathcal{T}$ compared to more complex reflections with respect to the boundary of $\mathcal{D}$ in \cite{AAFO}. Therefore, for small $j$ the adjusted statistic $\widehat{(s_0,\a)}_j$ is computationally simpler than the one in \cite{AAFO}.
\end{remark}

Now we are ready to formulate the main result.

\begin{theorem}\label{theo:ctl3}
Under the conditions of Theorem~{\rm \ref{theo:ctl2}}, the adjusted statistic $\widehat{(s_0,\a)}_j $ is a consistent asymptotically normal estimator of the parameter $(s_0,\a). $ When $j$ goes to $+\infty$, the  random vectors $\sqrt{m_j} \left( \widehat{(s_0,\a)}_j - (s_0,\a) \right)$
have the asymptotic bivariate centred Gaussian distribution  $\mathcal{N}(0, V_{s_0,\a} )$ with the covariance matrix $V_{s_0,\a}$ given by 
\begin{equation}\label{var_adj}
V_{s_0,\a}: = \frac{c \pi s_0^2  \int_{-c\pi}^{c\pi} \Big |\sum_{n\in\Z} \big | \wh{\psi} (\eta+2nc\pi)\big|^2 \Big |^2\,d\eta }
                  { 4\a^2 (1+2\ln s_0)^2 } 
             \begin{pmatrix} (V_{s_0,\a})_{11} & (V_{s_0,\a})_{12} \\ 
                             (V_{s_0,\a})_{12} & (V_{s_0,\a})_{22} \end{pmatrix} ,
\end{equation}
where
\begin{align}
(V_{s_0,\a})_{11}
& := \frac{(1-4\a\ln s_0)^2}{\left( \int_{\R}|\wh{\psi}(\eta)|^2\,d\eta \right)^2} 
  + \frac{8 s_0^4 (\ln s_0)^2}{\left( \int_{\R}\eta^2|\wh{\psi}(\eta)|^2\,d\eta \right)^2},
\nonumber \\
(V_{s_0,\a})_{12}
& := \frac{ (1-4\a\ln s_0)\a(4\a+2)s_0^{-1} }{ \left( \int_{\R}|\wh{\psi}(\eta)|^2\,d\eta \right)^2 } 
  - \frac{8\a s_0^3\ln s_0}{\left( \int_{\R}\eta^2|\wh{\psi}(\eta)|^2\,d\eta \right)^2},
\nonumber \\
(V_{s_0,\a})_{22}
& := \frac{ \a^2(4\a+2)^2 s_0^{-2} }{ \left( \int_{\R}|\wh{\psi}(\eta)|^2\,d\eta \right)^2 } 
  + \frac{8\a^2 s_0^2}{\left( \int_{\R}\eta^2|\wh{\psi}(\eta)|^2\,d\eta \right)^2}.\nonumber
\end{align}
\end{theorem}
\begin{proof1}{\it of Theorem~{\rm \ref{theo:ctl3}.}}\
The feasible region $\mathcal{D}$ is an open set.
Therefore, it follows from (\ref{theo:ctl3:eq_lfgn}) that, for any $\delta>0$ and for almost all $\omega \in\Omega,$ there is $J(\omega,\delta)$ large enough such that for $j\ge J$ the random vector
$
\left( \frac{ \ovl{\de}_j ^{(2,m_j)} }{ \int_{\R}|\wh{\psi}(\eta)|^2\,d\eta } ~,~
       \frac{ \De\ovl{\de}_{j+1}^{(2,M_j)} }{ 2\int_{\R}\eta^2|\wh{\psi}(\eta)|^2\,d\eta } \right)
$
 belongs to the $\delta$-neighbourhood of $\Phi(s_0,\a).$ Notice that  $1/m_j\to 0$ when $j\to +\infty.$ Hence, for almost all $\omega \in\Omega$ there is $J(\omega)$ large enough such that for $j\ge J$ the image under $\mathcal{T}(\cdot,{1}/{m_j})$ of the vector $
 \left( \frac{ \ovl{\de}_j ^{(2,m_j)} }{ \int_{\R}|\wh{\psi}(\eta)|^2\,d\eta } ~,~
        \frac{ \De\ovl{\de}_{j+1}^{(2,M_j)} }{ 2\int_{\R}\eta^2|\wh{\psi}(\eta)|^2\,d\eta } \right)
 $ equals to the vector itself.
 
 Thus, for $j \to +\infty$ 
\[
\sqrt{m_j} \left| \mathcal{T}\Bigg( \frac{ \ovl{\de}_j ^{(2,m_j)} }{ \int_{\R}|\wh{\psi}(\eta)|^2\,d\eta },
                                    \frac{ \De\ovl{\de}_{j+1}^{(2,M_j)} }{ 2\int_{\R}\eta^2|\wh{\psi}(\eta)|^2\,d\eta } ,
                                    \frac{1}{m_j} \Bigg)\right.\]
                                   \begin{equation}
                                   \label{theo:ctl3:eq_cvgps} 
                  \left.\quad \quad \quad - \Bigg( \frac{ \ovl{\de}_j ^{(2,m_j)} }{ \int_{\R}|\wh{\psi}(\eta)|^2\,d\eta } ,
                           \frac{ \De\ovl{\de}_{j+1}^{(2,M_j)} }{ 2\int_{\R}\eta^2|\wh{\psi}(\eta)|^2\,d\eta } \Bigg)  
           \right|
\mathop{\xrightarrow{\,a.s.\,}}0,
\end{equation}
where $|\cdot|$ is the Euclidean norm on $\R^2.$ Note that (\ref{theo:ctl3:eq_cvgps}) holds for any norm and any normalising factor, not only $\sqrt{m_j}$, because the difference almost surely vanishes for $j$ larger than some random $J.$

Hence, by (\ref{theo:ctl3:eq_lfgn}) and (\ref{theo:ctl3:eq_cvgps})
\[ 
\mathcal{T}\Bigg( \frac{ \ovl{\de}_j ^{(2,m_j)} }{ \int_{\R}|\wh{\psi}(\eta)|^2\,d\eta } ~,~
                         \frac{ \De\ovl{\de}_{j+1}^{(2,M_j)} }{ 2\int_{\R}\eta^2|\wh{\psi}(\eta)|^2\,d\eta } ~,~
                         \frac{1}{m_j} \Bigg)
\mathop{\xrightarrow{\,a.s.\,}} \Phi(s_0,\a),\ {j \to +\infty},
\]
which means that the vector $\mathcal{T}\Bigg( \frac{ \ovl{\de}_j ^{(2,m_j)} }{ \int_{\R}|\wh{\psi}(\eta)|^2\,d\eta } ~,~
                         \frac{ \De\ovl{\de}_{j+1}^{(2,M_j)} }{ 2\int_{\R}\eta^2|\wh{\psi}(\eta)|^2\,d\eta } ~,~
                         \frac{1}{m_j} \Bigg)$ is a consistent  estimator of $\Phi(s_0,\a).$
                  
Moreover, by multivariate Slutsky's lemma \cite[Theorem~2.7(iv)]{VanDerVaart}  it follows from (\ref{theo:ctl3:eq_cvgps})  and the central limit theorem (\ref{theo:ctl3:eq_clt1et2})  that for $j \to +\infty$  it holds
\begin{equation}
\label{theo:ctl3:eq_clt1et2_tronc} 
\sqrt{m_j} \left( \mathcal{T}\Bigg( \frac{ \ovl{\de}_j ^{(2,m_j)} }{ \int_{\R}|\wh{\psi}(\eta)|^2\,d\eta } ~,~
                                    \frac{ \De\ovl{\de}_{j+1}^{(2,M_j)} }{ 2\int_{\R}\eta^2|\wh{\psi}(\eta)|^2\,d\eta } ~,~
                                    \frac{1}{m_j} \Bigg)
                  - \Phi(s_0,\a) \right)  \mathop{\xrightarrow{\,d\,}}\mathcal{N}(0,V_{\V_1}),
\end{equation}
where 
\[V_{\V_1} := \V_1 \begin{pmatrix} \frac{1}{\left( \int_{\R}|\wh{\psi}(\eta)|^2\,d\eta \right)^2} & 0 \\ 
                                               0 & \frac{1}{2\left( \int_{\R}\eta^2|\wh{\psi}(\eta)|^2\,d\eta \right)^2} \end{pmatrix} .              
\]
%

The continuity of $\Phi^{-1}$ implies that the estimator $\widehat{(s_0,\a)}_j $ is consistent
\[ 
\widehat{(s_0,\a)}_j \mathop{\xrightarrow{\,a.s.\,}} (s_0,\a),\ j \to +\infty.
\]

As the central limit theorem in (\ref{theo:ctl3:eq_clt1et2_tronc}) can be rewritten as 
\[ 
\sqrt{m_j} \left( \Phi\big( \widehat{(s_0,\a)}_j \big) - \Phi(s_0,\a) \right)  
\mathop{\xrightarrow{~d~}} \mathcal{N}(0,V_{\V_1}),\ {j \to +\infty},
\]
then to obtain the asymptotic distribution of the estimator $\widehat{(s_0,\a)}_j$ around the parameter of interest $(s_0,\a)$ one can use the delta method with the inverse function $\Phi^{-1}$. 

To justify it one has to check that $\Phi^{-1}$ is differentiable at the point $\Phi(s_0,\a).$ By the inverse function theorem, the derivative $D(\Phi^{-1})(\Phi(s_0,\a))$ exists if the Jacobian $D\Phi$ of the function $\Phi(\cdot,\cdot)$ at the point $(s_0,\a)$ is  invertible. In this case it holds  $D(\Phi^{-1})(\Phi(s_0,\a))= \left( D\Phi(s_0,\a) \right)^{-1}.$

Notice that for any $(s_0,\a) \in (1,+\infty)\times(0,1/2)$  it holds
\begin{equation}\label{DPhi} 
D\Phi(s_0,\a) = 
                s_0^{-4\a-2} \begin{pmatrix} -4\a s_0             & -4 s_0^2 \ln s_0   \\ 
                                             \a (-4\a-2)s_0^{-1}  & 1 -4\a \ln s_0     \end{pmatrix} .
\end{equation}

Thus,  since $s_0>1,$
\[
\det\left( D\Phi(s_0,\a) \right) = s_0^{-8\a-4} \left( -4\a s_0 -8\a s_0 \ln s_0 \right)
                                 = -4\a s_0^{-8\a-3} (1+2\ln s_0) \neq 0 
\]
and the Jacobian matrix is invertible. 

Therefore, by the multivariate delta method (see, for example, \cite[Theorem~3.1]{VanDerVaart}) 
\[
\sqrt{m_j} \left( \widehat{(s_0,\a)}_j \big) - (s_0,\a) \right)  \mathop{\xrightarrow{\,d\,}}_{} \mathcal{N}(0, V_{s_0,\a} ),\ j \to +\infty,\]
where 
\begin{equation}\label{varV}
V_{s_0,\a} := \left( D\Phi(s_0,\a) \right)^{-1} V_{\V_1} \left( \left( D\Phi(s_0,\a) \right)^{-1} \right)^T.
\end{equation}

The covariance matrix given by (\ref{varV}) can be explicitly computed. It follows from (\ref{DPhi}) that \[
\left( D\Phi(s_0,\a) \right)^{-1}
=-\frac{s_0^{4\a+1}}{4\a (1+2\ln s_0)} \begin{pmatrix}  1 -4\a \ln s_0     & 4 s_0^2 \ln s_0  \\ 
                                                        \a (4\a+2)s_0^{-1} & -4\a s_0         \end{pmatrix}.
\]
Hence,
\[
V_{s_0,\a} = \frac{s_0^{8\a+2}\V_1}{16\a^2 (1+2\ln s_0)^2} 
             \begin{pmatrix} 1 -4\a \ln s_0     & 4 s_0^2 \ln s_0 \\ 
                             \a (4\a+2)s_0^{-1} & -4\a s_0        \end{pmatrix}    
\]
\[\times     \begin{pmatrix} \frac{1}{\left( \int_{\R}|\wh{\psi}(\eta)|^2\,d\eta \right)^2} & 0 \\ 
                             0 & \frac{1}{2\left( \int_{\R}\eta^2|\wh{\psi}(\eta)|^2\,d\eta \right)^2} \end{pmatrix}        \begin{pmatrix} 1-4\a\ln s_0   & \a (4\a+2)s_0^{-1} \\ 
                             4 s_0^2\ln s_0 & -4\a s_0        \end{pmatrix}.
\]
The straightforward matrix multiplication and application of (\ref{lem:sig-de2:eq2}) give (\ref{var_adj}), which completes the proof. 
\end{proof1}

\section{Numerical examples}\label{sec_6}

This section provides some numerical examples to illustrate and specify the general theoretical results from the previous sections.

The main theoretical results were obtained for general filter transforms and involve some complex functionals of the filters. The following two examples demonstrate that these results can be easily specialized for specific filters/wave\-lets and are  feasibly computable.\\

\begin{figure}[!htb]
\centering
\vspace{-0.3cm}
    \subfloat[Shannon father wavelet case]{\label{fig23a}
\includegraphics[trim={0.6cm 0cm 0cm 0cm},clip,width=0.45\textwidth, height=0.32\textheight]{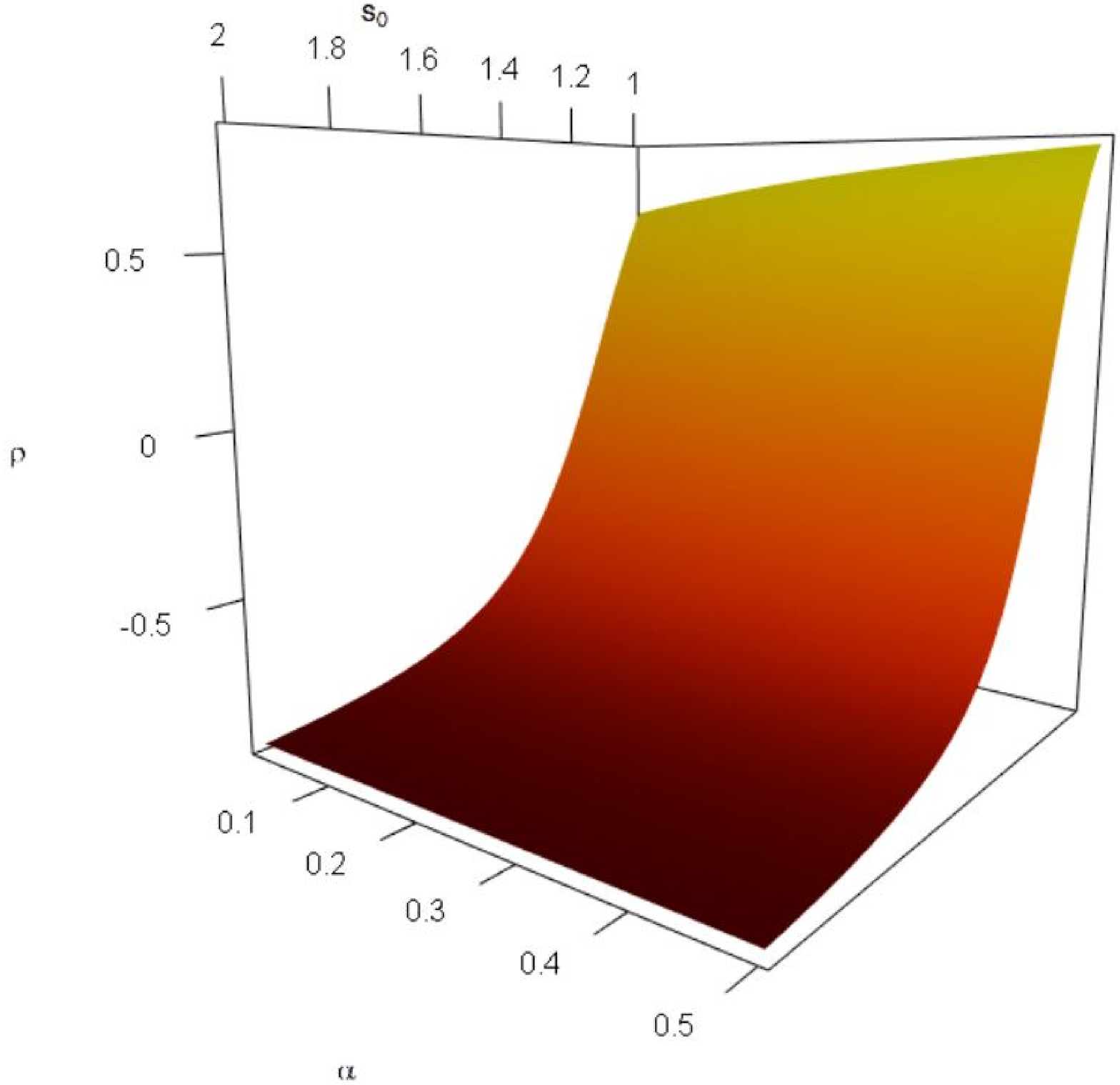}}
\centering
\vspace{-0.3cm}
    \subfloat[Meyer father wavelet case]{\label{fig23b}
\includegraphics[trim={0cm 0cm 1.5cm 0cm},clip,width=0.45\textwidth, height=0.32\textheight]{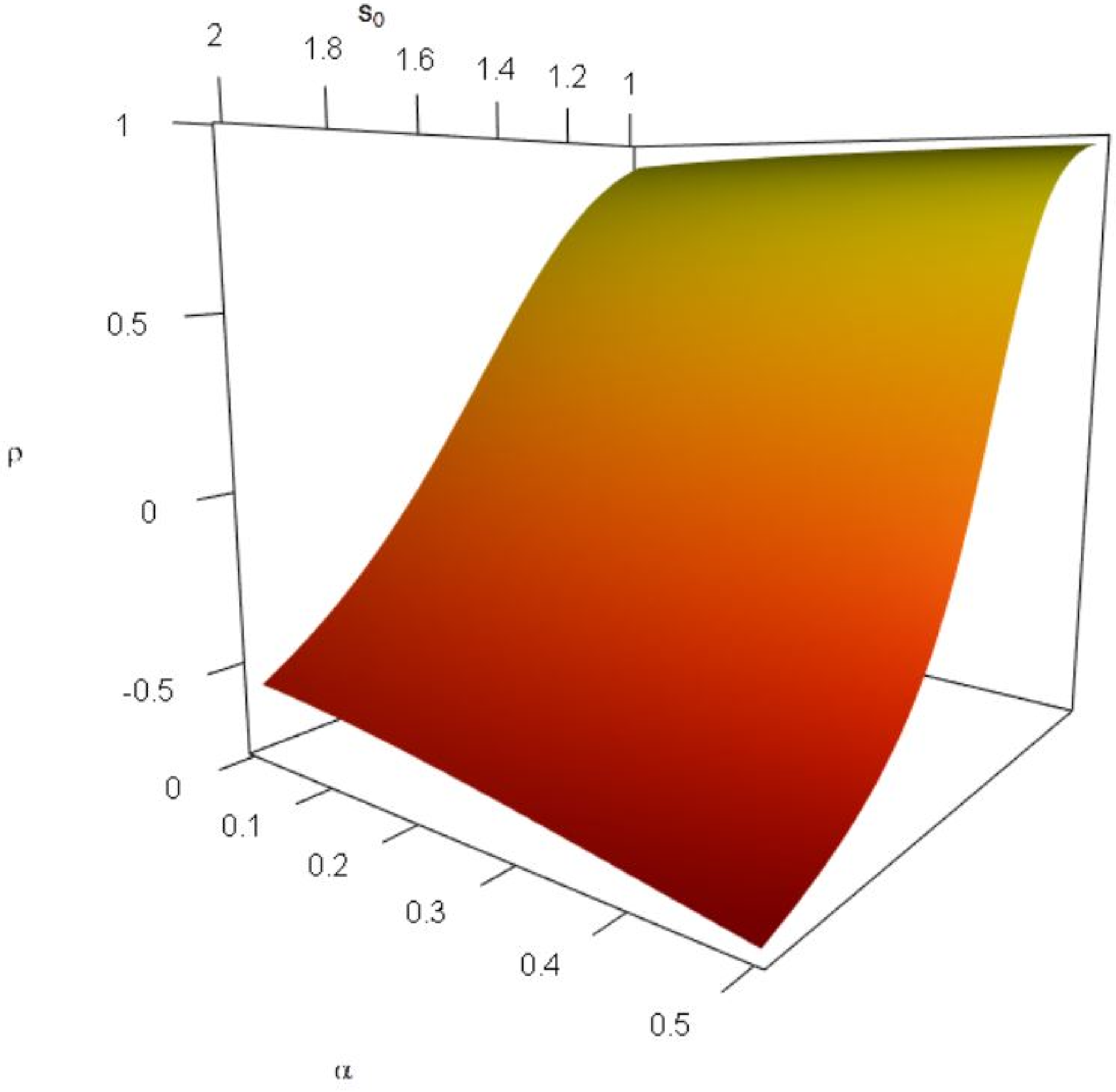}}
\vspace{0.4cm}
\caption{Asymptotic correlation of $\widehat{s_{0}}$ and $\widehat{\alpha}$.}
\end{figure}

\begin{Example}\label{ex1} Let us consider the Shannon father wavelet

\[ \psi_{f}(t) = {\rm sinc}(\pi t) :=
   {\begin{cases}\frac{\sin{(\pi t)}}{\pi t},& t \neq 0, \\
   1,& t=0.
   \end{cases}}
\]
Its Fourier transform is
\[ {\widehat{\psi}_{f}}{(\eta)} = I_{[-\pi, \pi]}{(\eta)} :=
   {\begin{cases}1,& \eta \in [-\pi, \pi], \\
   0,& \eta \notin [-\pi, \pi].
   \end{cases}}
\]
It is clear that Assumption~{\rm \ref{Assumption_2}} is satisfied.
The corresponding integrals are
\[\int_{\mathbb{R}}{{\left|{\widehat{\psi}_f}{(\eta)}\right|}^2 d\eta} = 2\pi \quad \quad \text{and} \quad \quad  \int_{\mathbb{R}}{{\eta}^2{\left|{\widehat{\psi}_f}{(\eta)}\right|}^2 d\eta} = {\frac{2}{3}}{\pi}^2.
\]
Let $I(c)$ denote the integral
\[ I(c) := \int_{-c\pi}^{c\pi}{\left|\sum_{n \in Z}{\left|{\widehat{\psi}_{f}}(\eta+2nc\pi)\right|}^2\right|}^2 d\eta = \int_{-c\pi}^{c\pi}{\left|\sum_{n \in Z} I_{[-\pi, \pi]}(\eta + 2nc\pi)\right|}^2d\eta.
\]

Then, for $c \geq 1$ one gets $I(c) = 2\pi.$

If $c < 1,$ by solving the inequality $c\pi + 2n^{*}c\pi \leq \pi$ we obtain 
$n^{*} =  \left[\frac{1-c}{2c} \right].$
Then, the solution of ${\eta}^{*} + 2(n^{*}+1)c\pi = \pi$ is ${\eta}^{*}={\pi\left(1-2c\left(1+\left[\frac{1-c}{2c}\right]\right)\right)}.$
Therefore, for ${\eta}^{*} < 0$ it holds
\[ I(c) = \int_{{\eta}^{*}}^{{-\eta}^{*}}(2n^{*}+1)^2d\eta +
          2\int_{{-c\pi}}^{{\eta}^{*}}(2n^{*}+2)^2d\eta\]
          \[ =
          -2{\eta}^{*}(2n^{*}+1)^2 + 2(c\pi + {\eta}^{*})(2n^{*}+2)^2 \]
and for ${\eta}^{*} \geq 0$
\[ I(c) = \int_{{-\eta}^{*}}^{{\eta}^{*}}(2n^{*}+3)^2d\eta +
          2\int_{{-c\pi}}^{{-\eta}^{*}}(2n^{*}+2)^2d\eta\]
          \[ =
          2{\eta}^{*}(2n^{*}+3)^2 + 2(c\pi - {\eta}^{*})(2n^{*}+2)^2. \]
Thus,
\[I(c) =
    \begin{cases}2\pi,& c \geq 1, \\
    2\left|{\eta}^{*}\right|(2n^{*}+2+{\rm sign}({\eta}^{*}))^2 + 2\left(c\pi - \left| {\eta}^{*}\right|\right)(2n^{*}+2)^2,& c <1.
    \end{cases}
\]
Hence, one can explicitly compute the covariance matrix $V_{s_{0},\alpha}$ in Theorem~{\rm \ref{theo:ctl3}}.
For example, the correlation of the components of the asymptotic vector equals 
\[ \rho = \frac{\frac{1}{4{\pi}^2s_{0}}(1-4\alpha\ln{s_{0}})\alpha(4\alpha+2)-\frac{18}{{\pi}^4}\alpha{s_{0}}^3\ln{s_0}}{\sqrt{{\Big(\frac{1}{4{\pi}^2}{(1-4\alpha\ln{s_{0}})}^2+{\frac{18}{{\pi}^4}}{s_{0}}^4{(\ln{s_0})}^2\Big)}{\Big(\frac{1}{4{\pi}^2{s_{0}}^{2}}{\alpha}^2(4\alpha+2)^2+{\frac{18}{{\pi}^4}}{\alpha}^2{s_{0}}^2\Big)}}}
\]
and is plotted in Figure~{\rm \ref{fig23a}} as a function of $s_{0}$ and $\alpha$. The plot shows that the components are highly correlated if $s_{0}$ is close to 1 and their correlation decreases as $s_{0}$ increases.
\end{Example}

\begin{Example}\label{ex2} Let us consider the Meyer father wavelet~{\rm \cite{Meyer:1992}}. It satisfies Assumption~{\rm \ref{Assumption_2}} as its Fourier transform equals
\[ {\widehat{\psi}_{f}}{(\eta)} = 
   {\begin{cases}1,& |\eta| \leq \frac{2\pi}{3}, \\
   \cos{\left(\frac{\pi}{2}\,\nu{\left(\frac{3|\eta|}{4\pi}-1\right)}\right)},& \frac{2\pi}{3} \leq |\eta| \leq \frac{4\pi}{3}, \\
   0,& otherwise,
   \end{cases}}
\]
where the function $\nu{(\cdot)}$ can be selected as
\[ \nu{(x)} = 
   {\begin{cases}0,& x < 0, \\
   x,& x \in [0, 1], \\
   1,& x > 1.
   \end{cases}}
\]
Its integrals are 
\begin{equation}\label{MeyerInt}
\int_{\mathbb{R}}{{\left|{\widehat{\psi}_f}{(\eta)}\right|}^2 d\eta} = 2\pi \quad  \quad \text{and} \quad \quad \int_{\mathbb{R}}{{\eta}^2{\left|{\widehat{\psi}_f}{(\eta)}\right|}^2 d\eta} = {\frac{8}{9}}{\pi}{({\pi}^2-2)}.    
\end{equation}
For example, for $c > \frac{4}{3}$ one can easily compute that
\[ I(c) = \int_{{-4\pi}/3}^{{4\pi}/3}{\left|{\widehat{\psi}}_{f}{(\eta)}\right|}^4 d\eta = {\frac{11}{6}}{\pi},
\]
which with {\rm (\ref{MeyerInt})} completely specifies the covariance matrix $V_{s_{0},\alpha}$. The corresponding correlation is shown in Figure~{\rm \ref{fig23b}} as a function of $s_{0}$ and $\alpha$.

Comparing it with Figure~{\rm \ref{fig23a}}, one can conclude that filters from Examples~{\rm\ref{ex1}} and~{\rm\ref{ex2}} produce similar correlation structures of the components of the asymptotic bivariate vector in Theorem~{\rm \ref{theo:ctl3}}. However, for the case of the Meyer father wavelet, the components exhibit higher correlations than for the Shannon one. 
\end{Example}

The following example continues simulation studies from \cite{AAFO}. Simulations in~\cite{AAFO} demonstrated consistency of the filter-based estimators of the cyclic and long-memory parameters. In Example~\ref{example3}, we examine their asymptotic normality. 

Note that the results in this paper were derived for functional time series with continuous time. For computer simulations, one has to use discretized processes on finite grids. In the available literature, it is usually assumed that the corresponding discretization error is negligible with respect to the estimation error. In many cases, it can be rigorously proven, see for example, \cite{Alodat:2020} and \cite{Ayaber:2011}.

\begin{Example}\label{example3} In this example the Mexican hat wavelet was used as a filter. This wavelet and its Fourier transform are  defined by, see {\rm \cite{Chun:2010}},  \[ \psi(t)=\frac{2}{\sqrt{3\sigma}\pi^\frac{1}{4} }\left( 1-\left(\frac{t}{\sigma}\right)^2\right) \mathrm{e}^{-\frac{t^2}{2\sigma^2}} \quad \mbox{\rm and}\quad \widehat \psi(\eta) =\frac{\sqrt{8} \pi^\frac{1}{4} \sigma^\frac{5}{2}} {\sqrt{3}} \eta^2 \mathrm{e}^{-\frac{\sigma^2 \eta^2}{2}}.\]
The value $\sigma=1$ was used for computations. The corresponding integrals are
\[\int_{\mathbb{R}}{{\left|{\widehat{\psi}}{(\eta)}\right|}^2 d\eta} = 2 \quad \quad \text{and} \quad \quad  \int_{\mathbb{R}}{{\eta}^2{\left|{\widehat{\psi}}{(\eta)}\right|}^2 d\eta} = 10.
\]
The Fourier transform $\widehat \psi(\eta)$ does not have a finite support, but has light tails that rapidly approaches zero  when $\eta\to +\infty.$ 

As  $X(t), t\in \mathbb{Z},$ we selected the Gegenbauer random process,  see {\rm \cite{Espejo:2015}}. This stochastic process is defined by the following difference equation 
\[
 \Delta^{d} _{u}X(t)=\varepsilon(t), \quad |u|\leq1, \ 0<|d|<1/2,
\]
where $\varepsilon(t)$ is a zero-mean white noise with the common variance $E (\varepsilon^2(t))= \sigma^2_{\varepsilon}.$

 The fractional difference operator $ \Delta^{d} _{u}$ is  given by 
\[\Delta^{d} _{u}= (1-2uB+B^2)^{d},\]
where $B$ denotes the time backward-shift operator, i.e. $B X(t)= X(t-1).$  

To simulate realizations of $X(t)$ we used truncated sums of the following infinite moving average representation of the Gegenbauer random process
\begin{equation}
 X(t)= \sum_{n=0}^{\infty}  C_{n}^{(d)} (u) \varepsilon(t-n),\quad  t\in \mathbb{Z},\label{GRF}\end{equation}
  with the coefficients given by the Gegenbauer polynomial 
\begin{equation*}
 C_{n}^{(d)} (u)= \sum_{k=0}^{[n/2]} (-1)^{k}\frac{(2u)^{n-2k}\Gamma (d-k+n)}{k!(n-2k)!\Gamma (d)},
\end{equation*} 
where $[n/2]$ is the integer part of $n/2,$ and  $\Gamma(\cdot)$ is the gamma function.
 
The chosen for simulations parameters values $d=0.1$ and $u=0.3$ correspond to $s_{0}$ and $\alpha$ inside of the admissible region $\mathcal{D}.$  The realizations of $X(t)$  were approximated by truncated sums with 100 terms in {\rm (\ref{GRF})}.  To compute the statistics  $\ovl{\de}_j ^{(2,m_j)}$ and $\De\ovl{\de}_{j+1}^{(2,M_j)}$  the values $a_j=j,$ $b_{jk}=k,$  $\gamma_j=1,$  and $ m_j=a_j^9,$ $j=1,...,7,$ were used.  It was shown in  {\rm \cite{AAFO}} that these values satisfy the assumptions of the theorems.

For $j=7,$ the subplots in Figures~{\rm\ref{fig24a}} and {\rm \ref{fig24b}} show Q-Q plots of the first two normalised statistics
\[S_1 := \sqrt{m_j}{\left(\frac{ \ovl{\de}_j ^{(2,m_j)} }{ \int_{\R}|\wh{\psi}(\eta)|^2\,d\eta} -
                         s_0^{-4\a}\right)} \] \text{and}
\[ S_2 :=  \sqrt{m_j}{\left(\frac{\De\ovl{\de}_{j+1}^{(2,M_j)}}{2\int_{\R}\eta^2|\wh{\psi}(\eta)|^2\,d\eta} - \a s_0^{-4\a - 2}\right)}.
\]
These plots demonstrate that these statistics have distributions close to Gaussian ones, which is also confirmed by the Shapiro-Wilk test for normality with the corresponding p-values 0.613 and 0.262. Moreover, the estimated correlation matrix
{\small $ \begin{pmatrix} 1  & 0.084 \\ 
                             0.084 & 1     \end{pmatrix}
$}
of these statistics and density ellipsoids in Figure~{\rm\ref{fig24c}} underpin the result in {\rm (\ref{theo:ctl3:eq_clt1et2})} about asymptotically bivariate normal distribution with uncorrelated components. Finally, Figure~{\rm\ref{fig24d}} gives density ellipsoids and realizations   of the random vector $\sqrt{m_{j}}\left( \widehat{(s_0,\a)}_j -(s_0, \a)\right)$ which suggest an asymptotically bivariate normal distribution as in Theorem~{\rm \ref{theo:ctl3}}.

\begin{figure}[!htb]
    \centering  \vspace{-0.3cm}
    \subfloat[Q-Q plot of $S_1$]{\label{fig24a}
    \includegraphics[trim={0cm 0cm 0cm 0cm},clip,width=0.45\textwidth, height=0.23\textheight]{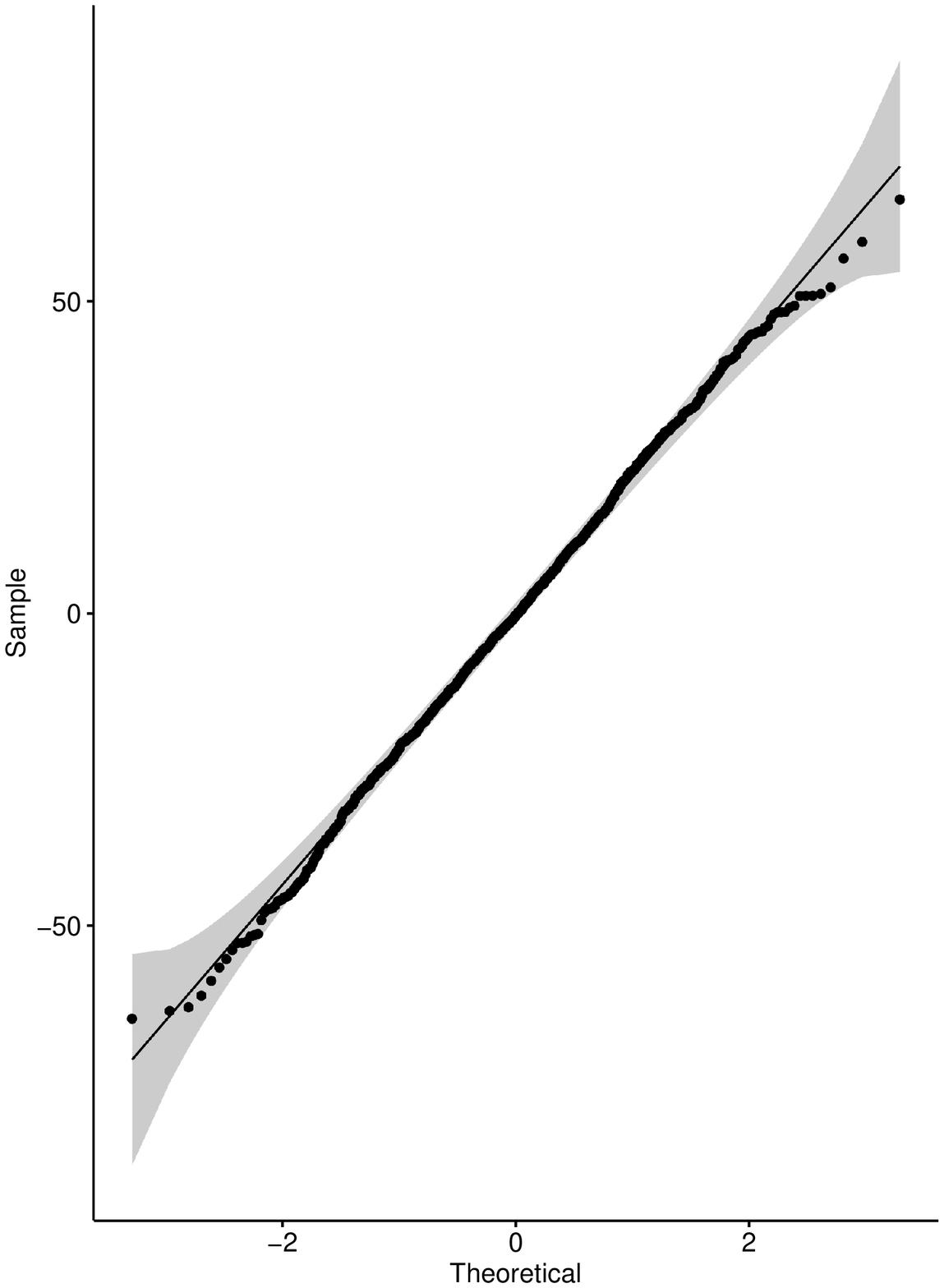}}
    \centering
    \subfloat[Q-Q plot of $S_2$]{\label{fig24b}
    \includegraphics[trim={0cm 0cm 0cm 0cm},clip,width=0.45\textwidth, height=0.23\textheight]{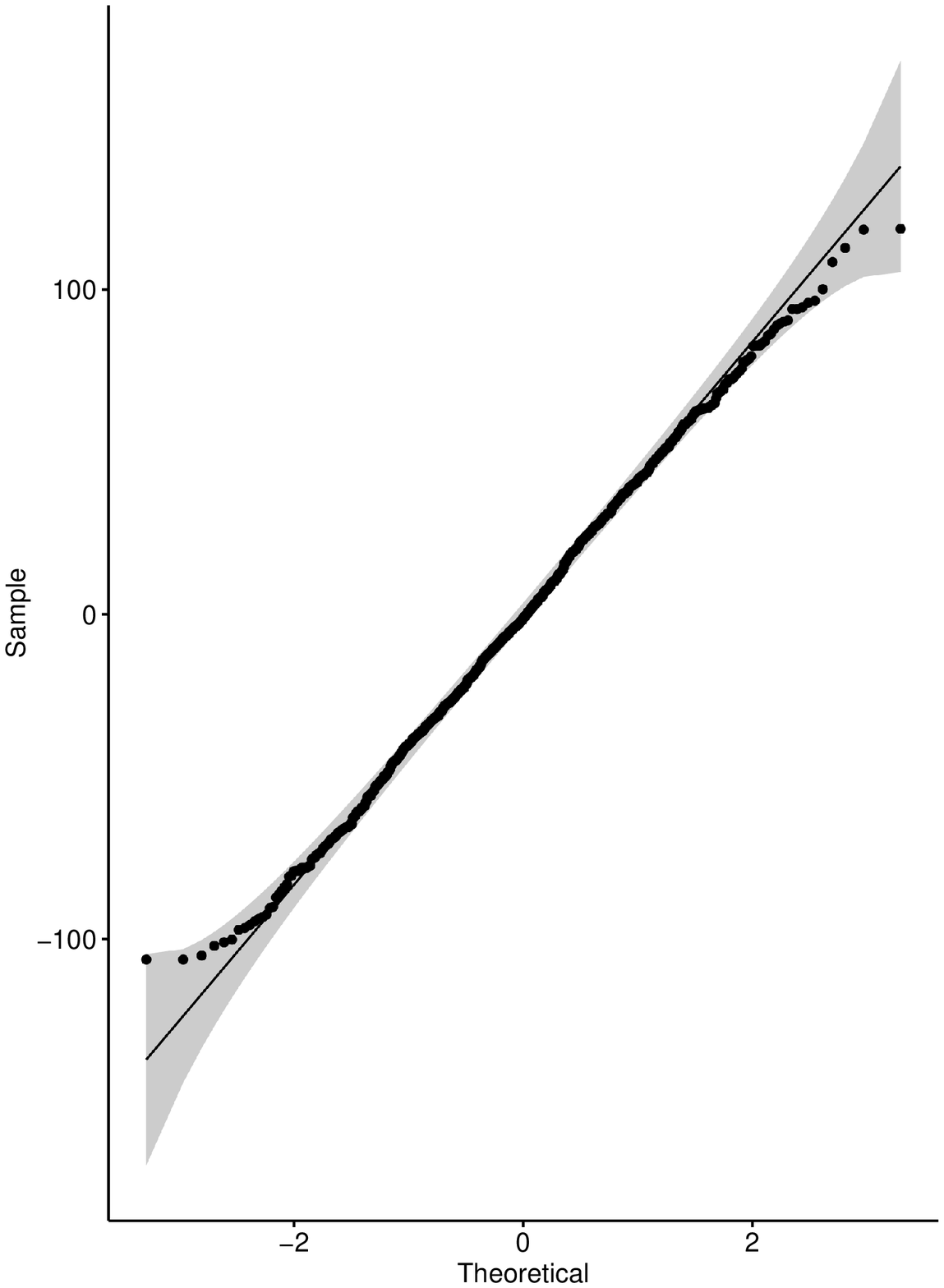}}\\
    \centering \vspace{-0.4cm}
    \subfloat[Density ellipsoid of ($S_1, S_2$)]{\label{fig24c}
    \includegraphics[trim={0.2cm 0cm 0.6cm 0cm},clip,width=0.45\textwidth, height=0.25\textheight]{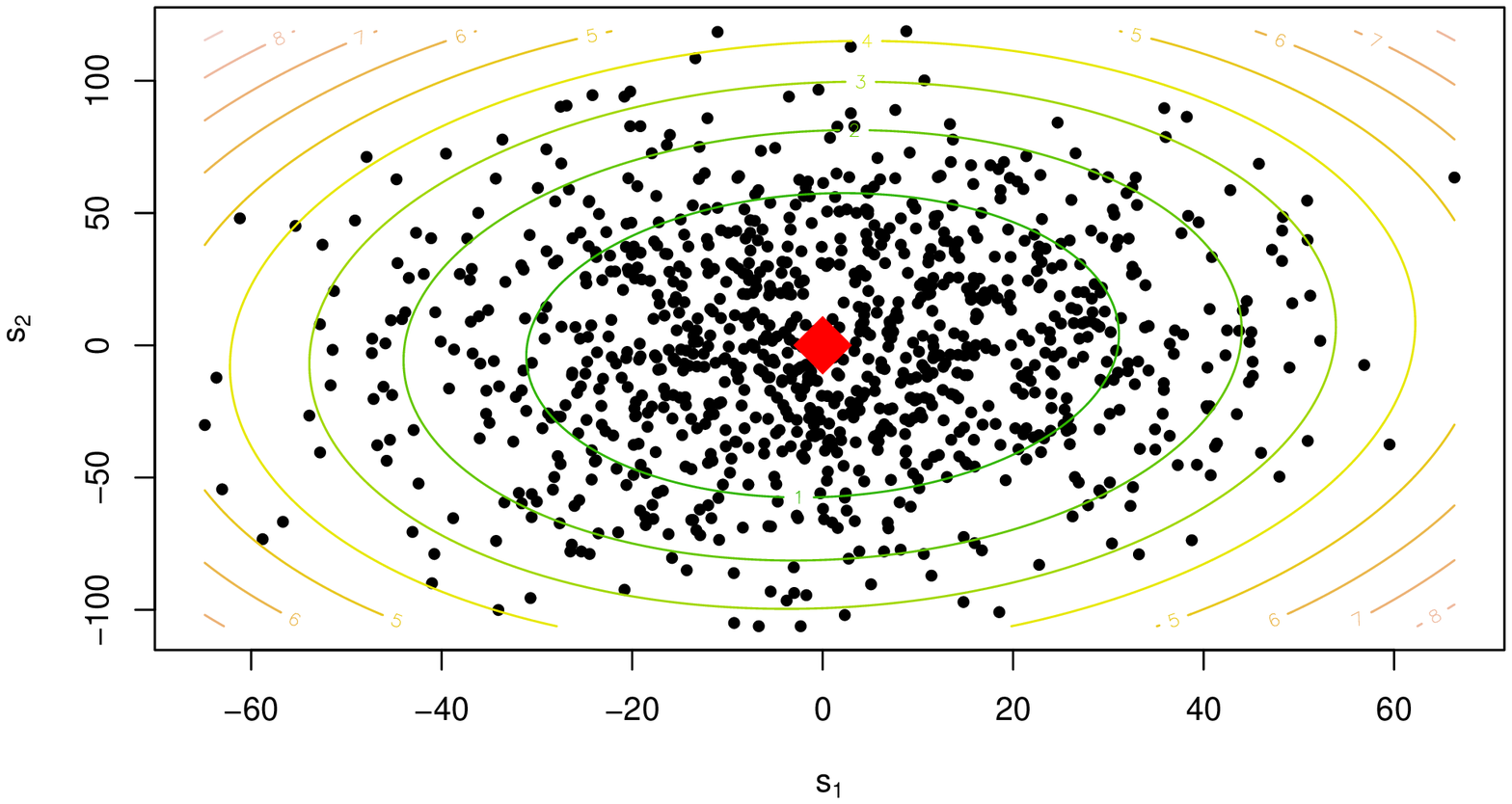}}
    \centering 
    \subfloat[\centering{Density ellipsoid of  \newline
    {${}\quad \sqrt{m_{j}}\left(\widehat{(s_0,\a)}_j -(s_0, \a)\right)$}}]{\label{fig24d}
    \includegraphics[trim={0cm 0cm 0.8cm 0cm},clip,width=0.45\textwidth, height=0.25\textheight]{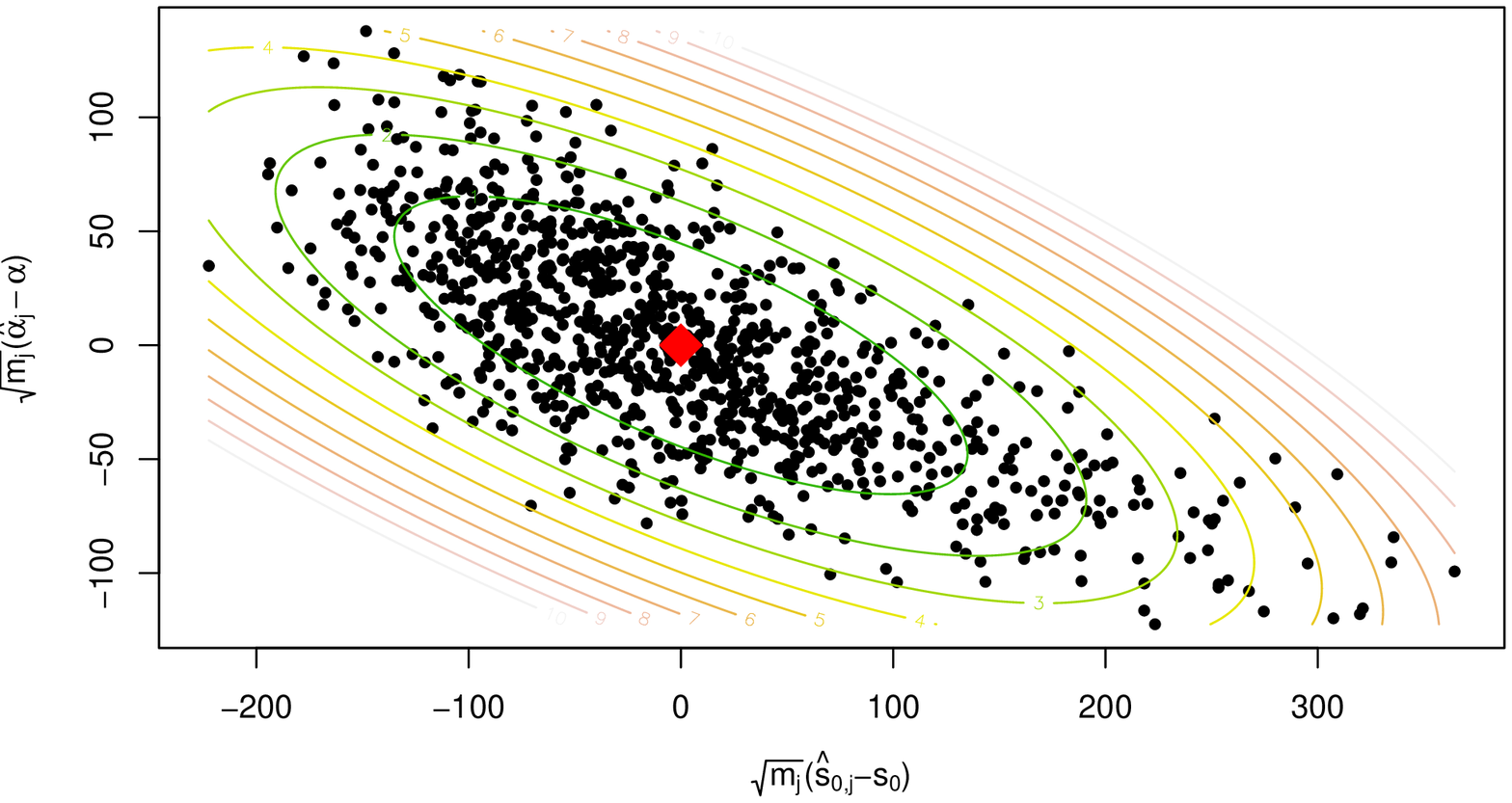}}
    \caption{Realizations of normalised statistics}
    \label{fig24}
\end{figure}

\end{Example}
The simulation studies suggest that the theoretical results are likely valid for wider classes of filters with light tails. They also demonstrate that the estimators exhibit approximately normal behaviour even for relatively small values of $j$. A separate publication will be devoted to comprehensive numerical studies.
\section{Conclusion}

The paper developed statistical inference of semiparametric models of functional time series. It was proved that the generalized filtered method-of-moment estimators of cyclic long-memory models are consistent and asymptotically normal. New adjusted simultaneous statistics were suggested and investigated. A rather general semiparametric class of models satisfies the assumptions of the theorems. In particular, Gegenbauer-type processes belong to this class.

Some interesting areas for future investigations are:

\begin{itemize}
    \item[--] Applying the approach to the case of multiple singularities, see~\cite{Arteche:2020, Klyolen:2012};
    \item[--] Adapting the methodology to models with other types of spectral singularities;
    \item[--] Investigating discretization errors for the case when $X(t)$ is observed on a finite grid, see~\cite{Ayaber:2011, Bardet:2010};
    \item[--] Investigating the case of random fields, i.e. when the index set of $X(t)$ is multidimensional, see~\cite{Ayache:2018, Espejo:2015, Klyolen:2012};
    \item[--] Continuing simulation studies to empirically compare the proposed approach with least squares and likelihood-type methods, see~\cite{RePEc:pra:mprapa:96313, Ferrara:2001, Whitcher:2004}.
\end{itemize}

\section*{Acknowledgements}
We are thankful to Professor Domenico Marinucci for attracting our attention to this research problem. 

Antoine Ayache is grateful to the Sydney Mathematical Research Institute at the University of Sydney for having financially supported his 7 weeks visit to Australia in~2019, and for the position of visiting researcher at the University of Sydney. Also, he is thankful to the School of Engineering and Mathematical Sciences at the La Trobe University for the honorary visiting professorship for one month at this university.

Andriy Olenko is grateful  to Laboratoire d'Excellence, Centre Europ\'{e}en pour les Math\'{e}matiques, la Physique et leurs interactions (CEMPI, ANR-11-LABX-0007-01), Laboratoire de Math\'{e}matiques Paul Painlev\'{e}, France, for support and giving him the opportunity to pursue research at the Universit\'{e} Lille for two months.

Ravindi~Nanayakkara and Andriy Olenko were partially supported under the Australian Research Council's Discovery Projects funding scheme (project DP160101366). 

This research includes computations using the Linux computational cluster Gadi of the National Computational Infrastructure (NCI), which is supported by the Australian Government and La Trobe University.

\bibliographystyle{imsart-number}

\bibliography{refs}

\end{document}